\documentclass{amsart}

\usepackage{geometry}
\usepackage[utf8]{inputenc}
\usepackage{amssymb}
\usepackage{amsmath}
\usepackage{amsfonts}
\usepackage{mdframed}
\usepackage{subcaption}
\usepackage{todonotes}
\usepackage{enumitem}
\usepackage{color}
\usepackage[algo2e,ruled, vlined, noend]{algorithm2e}

\usepackage{scalerel}
\usepackage{algorithm}

\newcommand{\sU}{\mathsf{U}}
\newcommand{\cF}{\mathcal{F}}

\newcommand{\cG}{\mathcal{G}}
\newcommand{\cE}{\mathcal{E}}
\newcommand{\bbR}{\mathbb R}
\newcommand{\bbP}{\mathbb P}

\newcommand{\mn}{\mathfrak q}
\newcommand{\mm}{\mathfrak p}

\newcommand{\cPnn}{{\mathcal P}_{n \times n}}
\newcommand{\cPn}{{\mathcal P}_{n}}
\newcommand{\one}{\pmb {1}}

\newcommand{\cSmn}{{\mathcal S}_{\mm,\mn}}
\newcommand{\cMn}{{\mathcal M}_n}
\newcommand{\cMnn}{{\mathcal M}_{n \times n}}

\definecolor{darkred}{rgb}{.7,0,0}

\definecolor{darkgreen}{rgb}{0,0.7,0}

\definecolor{darkblue}{rgb}{0,0,0.7}

\title{Inverse Optimal Transport}
\author{Andrew M. Stuart}
\address{California Institute of Technology, 1200 E. California Blvd, Pasadena, CA 91125}
\email{astuart@caltech.edu}
\author{Marie-Therese Wolfram}
\address{University of Warwick, Coventry CV4 7AL, UK and RICAM, Austrian Academy of Sciences, Altenbergerstr. 66, 4040 Linz, AT}
\email{m.wolfram@warwick.ac.uk}
 
\begin{document}
\maketitle

\begin{abstract}
Discrete optimal transportation problems arise in various contexts in 
engineering, the sciences and the social sciences. Often the underlying 
cost criterion is unknown, or only partly known, and the observed optimal 
solutions are corrupted by noise. In this paper we propose a systematic 
approach to infer unknown costs from  noisy observations of optimal 
transportation plans. The algorithm requires only the ability to solve
the forward optimal transport problem, which is a linear program, and to
generate random numbers. It has a Bayesian interpretation, and may also be
viewed as a form of stochastic optimization.

We illustrate the developed 
methodologies using the example of international migration flows. Reported 
migration flow data captures (noisily) the number of individuals moving 
from one country to another in a given period of time. It can be interpreted 
as a noisy observation of an optimal transportation map, with costs related 
to the geographical position of countries. We use a graph-based formulation
of the problem, with countries at the nodes of graphs and non-zero weighted
adjacencies only on edges between countries which share a border.  We use 
the proposed algorithm to estimate the weights, which represent cost of 
transition, and to quantify uncertainty in these weights.
  \end{abstract}

\section{Introduction}
\label{sec:I}

\subsection{Background}
\label{ssec:B}

There are many problems in engineering, the sciences and the social sciences in which
an input is transformed into output in an optimal way according to a 
cost criterion. We are interested in problems where the transformation 
from input to output is  known, and the objective is to infer the cost criterion which drives this
transformation. Our primary motivation is optimal transport (OT) problems in
which the transport plan is known but the cost is not. More generally linear
programs in which the solution is known, but the cost function and
constraints are to be determined, fall into the category of problems to
which the methodology introduced in this paper applies. We illustrate the type
of problem of interest by means of an example. \\

\noindent {\bf Example: International Migration.}  Quantifying migration flows between countries
is essential to understand contemporary migration flow patterns. Typically two types of migration statistics are collected -- flow
and stock data. Migration stock data states the number of foreign born
individuals present in a country at a given time and is usually based on population censuses.
Stock data is available for almost all countries in the world. Migration flow data captures the number  of migrants entering and leaving (inflow and outflow, respectively) a country over the course 
of a specific period, such as one year, see \cite{UN2017}. It is collected by most developed countries, but no international standards are defined. For example the time of residence after which a person counts  as an international migrant varies from country to country. Because of the different definitions and data collection methods, 
these statistics can be hard to compare. International agencies, such as the United Nations Statistics Division or the Statistical Office of the 
European Union (Eurostat),  publish annual migration flow estimates. 
These estimates are often based on Poisson or Bayesian linear regression. For more information about the estimation of migration flows using flow or stock statics we refer to \cite{Raymer2011, RWFSB2013, Abel2014,AR2019}. For the purposes
of this paper the main issue to appreciate is that migration data is available,
but should be viewed as noisy.

Flow data is typically presented in an origin-destination matrix, 
in which the $(i,j)^{\rm{th}}$ off-diagonal entry contains the number of people moving from country $i$ to 
country $j$ in a given period of time. This origin-destination data can be 
reported by both the sending (S) and the receiving (R) countries. Hence two migration flow tables 
are available, often desegregated by sex and age groups. Table \ref{t:migflows} shows harmonized data, which was pre-processed to improve comparability,   
reported by 6 European countries for the period 2002-2007. The numbers of the sending 
and receiving countries vary significantly. For example Germany reported that $136,927$ people 
immigrated from Poland, while Poland reported $14,417$ individuals who left for Germany. These very different numbers naturally raise the question of  the true migration flows.
In many settings it is natural to place greater weight on receiving data rather than departure data. But even this data is not subject to uniform standards, and therefore providing
reliable estimates and quantifying uncertainty is of great interest.

We interpret the reported origin-destination data maps (when appropriately normalized) as a noisy estimate of a transport plan arising from an OT problem
with unknown cost. It is then natural to try and infer the transportation
cost as it carries information about the migration process. $\diamond$

\begin{table}[h!]
  \centering
\begin{tabular}{r r|*{6}{c}}
\textbf{From} & & \textbf{To} & & & & \\
  & & CZ & DE & DK & LU & NL & PL \\
  \hline
 CZ & R&  0&	9,218&	262&	4&	511&	45\\
& S&	0&	    560&	24&	3&	81&	583	\\  
 DE & R &	1,362&	0&	4,001&	454&	9,182&	2,876\\
& S&	8,104&	0&	   3,095&	1,686&	9,293&	100,827\\
DK & R &	46&	2,687&	0&	11	&475&	34\\
& S&	179&	2,612&	0&	1,387&	602&	833\\
 LU & R &	2&	2,282&	162&	0&	161&	5\\
& S&	13&	911&	99&	0&	97&	23\\
 NL & R&	255&	13,681&	864&	27&	0	&163\\
& S&	298&	10,493&	533&	191&	0&	1,020\\
PL & R&	1,608&	136,927&	2,436&	19&	5,744&	0\\
& S&	63&	14,417&	111&	23&	577&	0\\
\hline
\hline
 \textbf{Tot:}&S &3,273&	164,795&	7,725&	515&	16,073&	3,123 \\
&R &8,657&	28,993&	3,862&	2,041&	10,650&	103,286\\
\end{tabular}
\caption{Harmonized migration flow statistics for the period 2002-2007; see \cite{BREW2010}. } \label{t:migflows}
\end{table}

The preceding example serves as motivation, and we will come back to it 
throughout this paper. However we reemphasize that the proposed identification methodologies
that we introduce in this paper can be used for general inverse OT 
and linear programming problems; further examples will serve to illustrate
this fact.

\subsection{Literature Review}
\label{ssec:LR}

\noindent Optimal transport originates with the French mathematician Gaspard
Monge who, in 1781, investigated the problem of finding the most
cost-effective way to move a pile of sand 
to fill a hole of the same volume. Kantorovich introduced the modern (relaxed) formulation of the problem, in which mass can be split, in 1942.
In more mathematical terms Kantorovich considered the following setup: given two positive
measures (of equal mass) and a cost function, find the transportation map that
moves one measure to the other minimizing the transport cost. The
corresponding infimum induces a distance between these two measures -- the
so-called Wasserstein distance. The Wasserstein distance plays an important
role in probability theory, partial differential equations (PDE) 
and many other fields in applied mathematics \cite{V2013,S2015}. Furthermore 
the techniques and methodologies developed in OT have found application in
a variety of  scientific disciplines including data science, economics,
imaging and meteorology  \cite{GS2015}.

With the spread and application of OT
into different scientific disciplines the interest in computational
methodologies has increased. Commonly used numerical methods  broadly speaking
fall into two categories: linear programming \cite{dantzig2016linear}
and methods specific to the structure of OT. Linear programs are classic
problems which have been extensively studied in the field of optimization and
operations research. Many computational methodologies have been developed, such as the famous simplex
algorithm (and its many variants), the Hungarian algorithm and the auction
algorithm. All these methods work well for small to medium sized
problems, but are too slow in modern applications such as imaging or  supply
chain management. Recently a significant speed up, of linear programming,
was achieved by considering a regularized OT problem, leading to the
Sinkhorn algorithm (or variants thereof) in which an additional entropic 
regularization term is added to the objective function; this allows efficient
computation of the corresponding minimizer and induces a trade-off
between fidelity to the original problem, and computational speed.
This family of efficient algorithms resulted in the rapid advancement of 
computational OT in recent years, especially in the context of imaging and data science; see \cite{Peyre2017,Cuturi2013,reich2018data}.

Inverse problems for linear programming received considerable interest in the engineering literature. The
paper \cite{AO2001}, building on earlier work in \cite{zhang1996calculating},
studies the problem by seeking a cost function nearest to a 
given one in $\ell^p$ for which the given solution is an optimal linear program;
this problem is itself a linear program in the case $p=1.$
The formulation of an inverse problem for linear programming in \cite{DL2006}
took a slightly more general perspective, as it does not assume that the given
data necessarily arises as the solution of a linear program, and rather seeks
to minimize the distance to the solution set of a linear program. Recent
application of the inverse problem for linear programming may be found in \cite{SM2018},
for example.
These papers on inverse linear programming are foundational and have
opened up a great deal of subsequent research. However the methods in them do not
account in a systematic way for noise in the data provided, and
for the incorporation of prior information.
We address these issues by adopting a Bayesian formulation of the
inverse problem for linear programming, concentrating on OT 
in particular; the ideas are readily generalized to inverse linear programming
in general. The Bayesian approach not only allows for the quantification
of uncertainty, but also leads to new (stochastic) optimization methods.
An overview of the computational state of the art for Bayesian inversion may be found in
\cite{kaipio2006statistical}. The specific methods that we introduce
have the desirable feature that they require only solution of the forward
OT problem and the ability to generate random numbers.

\subsection{Our Contribution}
\label{ssec:C}

Our contributions to the subject of inverse problems within linear programming
are as follows.

\begin{itemize}

\item We formulate inverse OT problems in a Bayesian framework.

\item We provide a computational framework for solving inverse OT problems
in an efficient fashion.

\item We introduce graph-based cost functions for OT, using graph-shortest
paths in an integral way.

\item Graph-based OT has considerable potential for application,
and we introduce a new way of studying migration flow data using 
inverse OT in the graph-based setting.

\end{itemize}

We emphasize that, whilst the graph-based formulation of cost corresponds 
to a rather specific way of designing cost functions for discrete linear 
programs, the framework and algorithms developed in this paper apply
quite generally to inverse linear programming, and hence to OT in general. We develop the methodology in general,
using graph-based migration flow as a primary illustrative example.
In section \ref{sec:IOT} we define OT as a linear program, describe
the cost criteria considered, and formulate inverse OT in a Bayesian
setting. Section \ref{sec:A} presents algorithms for the forward and inverse OT 
problem and section \ref{sec:NR} contains numerical results.\\

We will use the following notation throughout this manuscript.
Let $|\cdot|$ and $\langle \cdot, \cdot \rangle$ denote the
Euclidean norm and inner-product on $\bbR^n$ and the Frobenius norm and
inner-product on $\bbR^{n \times n}.$ The spaces of probability matrices,
probability vectors and probability matrices with specified marginals 
are defined as
\begin{align*}
  &\cPnn=\Bigl\{B \in \bbR^{n \times n}: B_{ij} \ge 0, \sum_{i,j=1}^n B_{ij}=1\Bigr\},~~\cPn=\Bigl\{u \in \bbR^{n}: u_{j} \ge 0, \sum_{j=1}^n u_{j}=1\Bigr\},\\
  &\cSmn=\Bigl\{B \in \cPnn: B\one=\mm, B^T\one=\mn\,\,{\rm for}\,\, \mm,\mn \in \cPn\Bigr\} \text{ where }\one=(1,\cdots, 1)^T \in \bbR^n. 
  \end{align*}

\section{Inverse Optimal Transport} \label{sec:IOT}
In this section we introduce the forward OT problem and discuss specific cost criteria, before formulating the respective inverse OT problem in the Bayesian framework.

\subsection{Forward Problem}
\label{ssec:FP}
We consider two discrete probability vectors $\mn \in \cPn$ and $\mm \in \cPn$ and a given cost $C \in \cPnn$. Then the optimal transport problem corresponds to finding a  map
transporting $\mm$ to $\mn$ at minimal cost. Note that in OT the cost matrix
has non-negative entries, which can be normalized to be an element of $\cPnn$ without loss of generality.
The respective forward OT problem is to find
\begin{align}\label{e:ot}
T^* \in {\rm argmin}_{T \in \cSmn} \langle C,T \rangle.
\end{align}
Problem \eqref{e:ot} falls into the more general class of linear
programs. Linear programs (and their many variants) arise in
various specific settings -- such as
the earth mover's distance (EMD)\cite{Rubner2000} or cost network flows \cite{BPPH2011} 
-- in different scientific communities. The problem 
\eqref{e:ot} has, by virtue of being a specific class of linear programs, 
at least one solution; this solution lies on the boundary of the feasible set of solutions (defined by the equality constraints).
If the solution is unique then we define mapping 
$\cF: \cPn \times \cPn \times \cPnn \rightarrow \cPnn$ by
\begin{equation}
\label{eq:forward}
T^*=\cF(\mm,\mn,C).
\end{equation}
In the non-unique setting we define $\cF(\mm,\mn,C)$ to be a unique element 
determined by running a specific non-random algorithm for the linear program
to termination, started at a specific initial guess.

We now consider \eqref{e:ot} regularized by the addition of the discrete entropy,
an approach popularized in \cite{Cuturi2013,Peyre2017} and which has led to considerable
analytical and computational developments. The resulting problem is
\begin{align}\label{e:entropy}
H(T) = -\langle T, \log(T)\rangle+{\rm Tr}(T)=-\sum_{i,j=1}^n {T}_{i,j}(\log {T}_{i,j} - 1),
\end{align}
where the matrix logarithm operation is applied
elementwise. Then
\begin{align}\label{e:regot}
T^*_{\epsilon}={\rm argmin}_{T \in \cSmn} \Bigl(\langle C,T \rangle+\epsilon H(T)\Bigr).
\end{align}
This problem has a unique minimizer $T_{\epsilon}^*$, since $H(T)$ is strongly
convex. Following our previous notation we define the corresponding mapping by $\mathcal{F}_{\epsilon}: \cPn \times \cPn \times \cPn \rightarrow \cPnn$
\begin{align}
T^*_{\epsilon} = \cF_{\epsilon}(\mm,\mn,C).
  \end{align}
It is, in contrast to the optimal solution of \eqref{e:ot}, not sparse. It is known that solutions to \eqref{e:regot} converge to minimisers of \eqref{e:ot} as $\epsilon \rightarrow 0$.
Determining the rate of convergence is still an open problem.
The special structure of this regularized problem can be used to construct efficient splitting algorithms. These methods are based on the equivalent formulation of finding the projection of the joint coupling with respect to the Kullback-Leibler divergence
\begin{align*}
D_{\scaleto{KL}{3pt}}(T\|K) := \langle T, \log(T/K)\rangle-{\rm Tr}(T)+{\rm Tr}(K)=\sum_{i,j=1}^n T_{i,j}\log\frac{T_{i,j}}{K_{i,j}} - T_{i,j} + K_{i,j},
  \end{align*}
where the matrix logarithm and division operations are applied
elementwise and $K$ is the Gibb's kernel
\begin{align}\label{e:gibbs}
 K_{i,j}= \exp^{-\frac{C_{i,j}}{\epsilon}}.
\end{align}
In particular
\begin{align}\label{e:kldiv}
  T^*_{\epsilon} = {\rm argmin}_{T \in \cSmn} D_{\scaleto{KL}{3pt}}(T\|K).
\end{align}
The Kullback-Leibler divergence can be computed extremely efficiently using
proximal methods, yielding for example the celebrated Sinkhorn algorithm. We will
briefly outline the underlying ideas in Section \ref{ssec:cot}.

\subsection{Cost Criteria}
\label{ssec:GC}
\begin{figure}
  \centering
  \includegraphics[width=0.5\textwidth]{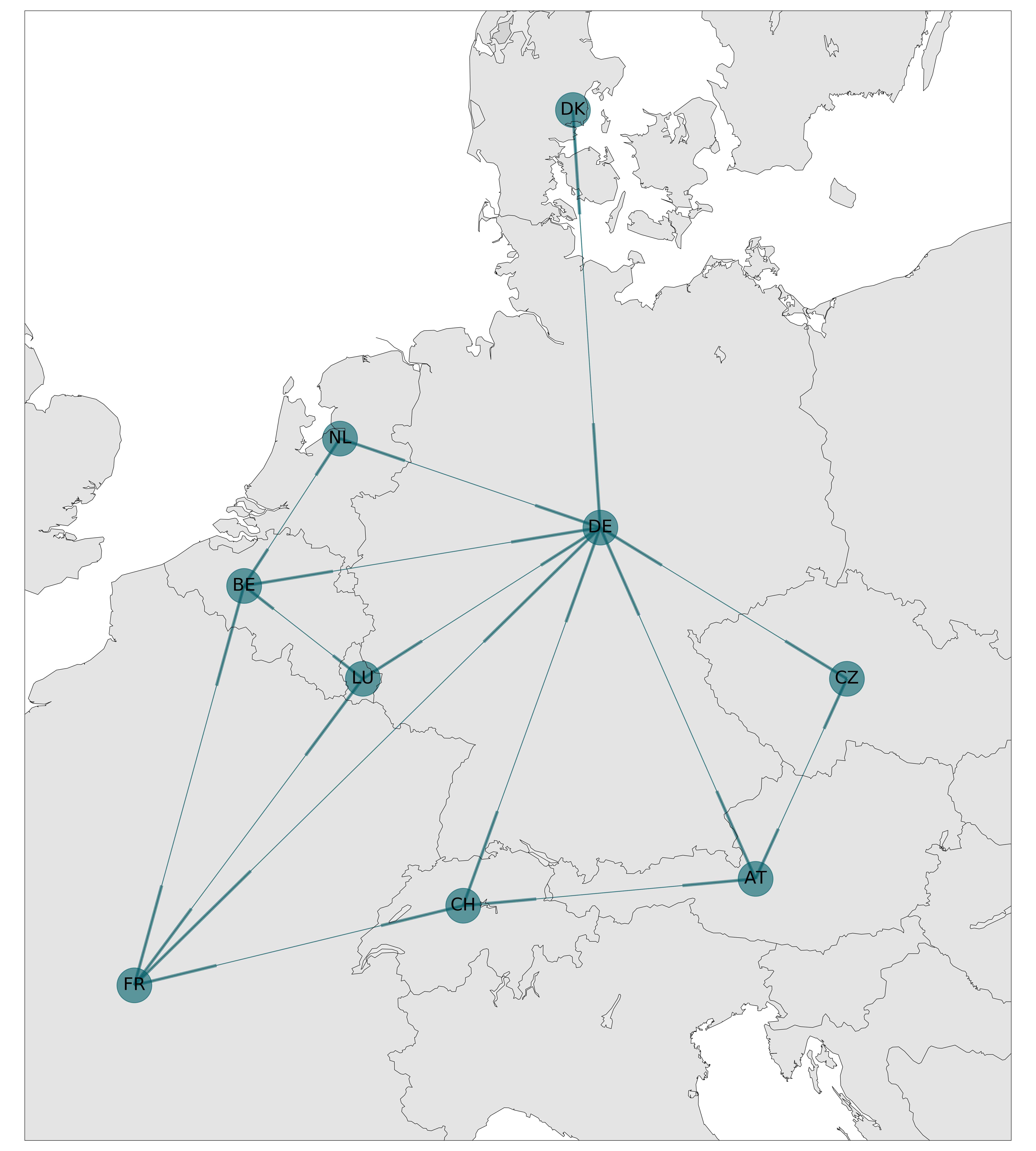}
  \caption{Network defined by the European countries used in our example.}
  \label{f:europe}
\end{figure}

Problems \eqref{e:ot} or \eqref{e:regot} are formulated for general cost matrices $C$ - the specific structure of $C$ depends on the application considered. We will investigate the behavior of the proposed methodologies for $C$ being:
\begin{enumerate}[label=(\roman*)]
 \item Toeplitz; \label{i:toeplitz}
 \item non-symmetric;\label{i:general}
 \item determined by an underlying graph structure.\label{i:graph}
\end{enumerate}
We assume that all individuals move, hence $T_{ii} = 0$ for all $i = 1, \ldots
n$ in all three cases, Therefore 'staying' is penalized by setting
\begin{align}\label{e:pen}
  C_{ii} = \bar{C} \gg 1\quad \text{ for all } i=1, \cdots n.
\end{align}
If $C$ is Toeplitz the cost depends on the difference between indices and $C$
has  $2n-3$ degrees of freedom.  Case \ref{i:general} corresponds to general
non-symmetric transportation cost, which in the context of migration flows
could include factors such as sharing the same language, the ratio of the
gross national income per capita or their EU membership. In case \ref{i:graph}
we assume that costs are related to an underlying discrete
structure. In the context of migration flows the geographical position of
countries defines an underlying graph with edges only between
countries which share a border; see Figure \ref{f:europe}.
We assume that the total transportation cost corresponds to the sum of the individual costs of moving from one country to another along edges of the
graph. In defining cost this way we are implicitly assuming that, between
the European countries studied here, migration is primarily via land.
This resulting discrete underlying structure, which relates the cost matrix to 
a directed graph representing the migration network between countries, is
detailed in the following.

Let $(V,E)$ be a directed graph with $n=|V|$ vertices and a (possibly non-symmetric) weighted
adjacency matrix $A \in \bbR^{n \times n}.$ We can then define a cost matrix $W \in \bbR^{n \times n}$ whose $(i,j)^{\rm th}$ entry $W_{i,j}$ is the shortest path cost of moving from vertex $i$ to $j$
according to the weighted adjacency matrix $A.$  Let $m$ be the number of non-zero
entries of $A$ and $f \in \bbR^m$ the vector defining the non-zero entries. 
Then we may define a mapping  $\cE$
such that $W=\cE(f).$ This $W \in \bbR^{n \times n}$ can then be normalized to
give a $C \in \cPnn$ and we may define the solution of the resulting OT 
problem via \eqref{eq:forward}.
For this graph-based cost the solution of the OT problem may be viewed as a function
of $\mm,\mn$ and $f$.
The minimal cost of moving between vertices of a graph can be computed using Dijkstra's algorithm, recalled in Section \ref{ssec:cot} below.

We define a similar mapping in the case of Toeplitz cost. Here the respective 
cost matrix $C$ has $2n-2$ free entries, before 
normalization to a probability vector and recalling that we fix the diagonal
to penalize not moving, 
and so we define a mapping $\mathcal{E}: f \in \mathbb{R}^{2n-2} \rightarrow 
\bbR^{n \times n}_+$; normalization then gives $C = \cMnn(\mathcal{E}(f))$.

\subsection{Inverse Problem}
\label{ssec:IP}
The inverse OT problem is to find $\mm,\mn$ and $C$ from the solution $T$
to the OT problem \eqref{e:ot}, or its regularized counterpart \eqref{e:regot}.
We tackle this problem by introducing a space of componentwise
positive and real-valued latent 
variables $u,v,W$ or $u,v,f$ which map to the unknowns $\mm,\mn \in \cPn$ 
and $C \in \cPnn$. It
is easier, and more natural, to specify priors in terms of these real-valued
latent variables.
To this end we introduce mappings from $\bbR^n_+$ into $\cPn$ and
from $\bbR^{n \times n}_+$ into $\cPnn$ as follows:
$\cMn: \bbR^n_+ \mapsto \cPn$ is defined by
$$\cMn(u)_j=u_j/(\sum_{\ell=1}^n u_{\ell})$$
and
$\cMnn: \bbR^{n\times n}_+ \mapsto \cPnn$ is defined by
$$\cMnn(W)_{i,j}=W_{i,j}/(\sum_{k,\ell=1}^n W_{k,\ell}).$$
Note that $\cMn(\lambda u)=\cMn(u)$ for all $\lambda \in \mathbb{R}$;
the same holds for $\cMnn.$
Then the forward problem \eqref{eq:forward} can be written as 
\begin{align}\label{e:iot}
  T^*=\cG(u,v,W):=\cF\bigl(\cMn(u),\cMn(v),\cMnn(W))
  \end{align}
or, in the case of graph-based cost or Toepliz cost, we have
\begin{align}\label{e:iot_graph}
  T^*=\cG(u,v,f):=\cF\bigl(\cMn(u),\cMn(v),\cMnn(\cE(f))).
\end{align}
This is readily generalized to the use of regularized optimal
transport as the forward model, simply replacing $\cF$ by
$\cF_{\epsilon}.$

We wish to invert the map $\cG$, given noisy observations of $T^*.$ Such problems are in general ill-posed, hence  suitably regularized versions have to be considered. Different approaches can be found in the literature -- we focus on the Bayesian framework, which allows us to estimate the posterior distribution of $u, v$ and $W$ (or $f$).

Depending on the structure of the cost matrix the inverse problem related to \eqref{e:iot} or \eqref{e:iot_graph} can be over- or underdetermined. 
We recall that in case of Toeplitz cost the matrix $C$ has $2n-3$ degrees of freedom. Then we have  $n^2-1$ equations for $4n-5$ unknowns (taking into account the normalization
of $u$, $v$ and $W$). Hence the inverse problem is overdetermined for $n>2$. If $C$ is a  general cost matrix with a set penalty on the diagonal, that is case \ref{i:general}, the cost matrix has
$n^2 - n$ degrees of freedom. In total we have  $n^2 + n -3$ unknowns, and
therefore the problem is underdetermined for $n > 2$.  For graph-based cost
(case \ref{i:graph}) the matrix $C$ has $m$ degrees of freedom and therefore the problem is underdetermined if
\begin{align}
2n+m-3 > n^2 - 1.
  \end{align}

\subsubsection{Likelihood}
We define a Bayesian formulation of the inverse problem, working in the
case where $u,v,W$ are the unknowns; the extension to $u,v,f$ as unknowns
is similar. We assume that the observed transport maps are corrupted by noise, in particular
\begin{align}\label{e:noise}
  T=\cG(u,v,W)+\eta
  \end{align}
where $\eta$ is a mean zero Gaussian random matrix with i.i.d. entries
of variance $\sigma^2.$
 In particular we want to find the conditional probability distribution of $ (u,v,W)$ given noisy observations $T^*$, that is $(u,v,W) \mid T^*$.
 The estimation is based on the model-data misfit function
\begin{align}\label{e:misfit}
\Phi(u,v,W;T) = \frac{1}{2\sigma^2}\lvert T - \cG(u,v,W))\rvert^2.
  \end{align}
Using Bayes' formula 
\begin{align}\label{e:bayes}
\bbP(u,v,W \mid T) = \frac{1}{\bbP(T)} \bbP(T \mid u,v,W) \bbP(u,v,W)
\end{align}
the probability of observing $T$ given a realization of $u, v$ and $W$, which is exactly the posterior distribution of $u, v$ and $W$, is given by
\begin{align}
  \bbP(T|u,v,W) \propto \exp\Bigl(-\frac{1}{2\sigma^2}|T-\cG(u,v,W)|^2\Bigr) = \exp\Bigl(-\Phi(u,v,W;T)\Bigr).
  \end{align}
In \eqref{e:bayes} $\bbP(u,v,W)$ corresponds to the prior information about $u,v$ and $W$. We assume that $u,v$ and $W$ have i.i.d. entries uniformly distributed
in $[0,1]$ and denote the set of vectors and matrices which satisfy this
componentwise constraint by $\sU.$ In view of the scale invariance of $\cMn(\cdot)$ 
the choice of unit interval $[0,1]$ is immaterial; any bounded interval $[0,\lambda]$ 
would deliver identical posterior on $u,v,W$.\\
Then the posterior distribution of $u,v$ and $W$ is given by
\begin{align}\label{e:posterior}
  \bbP(u,v,W|T)&=\frac{1}{Z} \exp\Bigl(-\frac{1}{2\sigma^2}|T-\cG(u,v,W)|^2 \Bigr)\one_{\sU}(u,v,W),
  \end{align}
with a normalization constant 
$$Z=\int_{\sU}\exp\Bigl(-\frac{1}{2\sigma^2}|T-\cG(u,v,W)|^2 \Bigr)du\, dv\, dW\,.$$

\noindent We can either sample from the posterior \eqref{e:posterior} (which corresponds to the full Bayesian approach) or maximize the
posterior probability \eqref{e:posterior}, which
leads to the optimization problem of minimizing $\Phi(u,v,W;T)$ over $\sU$. The first approach allows us to quantify uncertainty in the estimates of $u, v$ and $W$, the latter gives a single estimate. We discuss how to sample from the posterior, using a Random walk Metropolis (RwM) method in  Section \ref{ssec:MCMC}

\section{Algorithms For Inversion}
\label{sec:A}
\noindent In the following we present the numerical methods used in the computational
experiments in Section \ref{sec:NR}. Since the proposed Bayesian framework
requires the solution of an OT problem \eqref{e:ot} (or its regularized
version \eqref{e:regot}) in every iteration of the sampling algorithm,                                             
computational efficiency is essential. We start by presenting the solvers for the forward OT problem followed by the
Markov-Chain-Monte-Carlo methods used to sample from the posterior.

\subsection{Computational Optimal Transport}
\label{ssec:cot}

Numerical methods for linear programming go back to the seminal works of
Dantzig on the simplex method, see \cite{dantzig2016linear}.  Solutions to the linear program \eqref{e:ot} lie on the boundary of the feasible polytope, which is
defined by the constraints. The simplex method iterates over the vertices of this polytope to find the optimal solution, see \cite{NW2006}. The method works well in practice, however
examples in which the performance scales exponentially with the dimension of the problem, can be constructed. Different approaches to speed up computations have been proposed: for example
network simplex algorithms are based on the fact that specific linear programs can be formulated as minimization problems on graphs. The particular structure of the underlying graph can be
used to speed up the simplex method significantly. Further information on computational methods for linear programming can be found in \cite{BREW2010}. \\

\noindent More recently computational techniques, which are based on the regularized OT problem
\eqref{e:regot} have been proposed in the literature. These methods are
extremely efficient, since they are based on the formulation of the OT problem
in terms of the Kullback Leibler divergence \eqref{e:kldiv}. Its minimiser is given by
\begin{align*}
 T_{i,j} = a_i K_{i,j} b_j.
\end{align*}
Here $K$ is the Gibb's kernel \eqref{e:gibbs} and the vectors $a$ and $b$ satisfy the mass constraint
\begin{align}
 {\rm diag}(a) K \one = \mm \text{ and } {\rm diag }(b) K^T \one    = \mn.
\end{align}
This mass constraint can be enforced iteratively via
\begin{align}\label{e:sinkhorn}
 a^{(l+1)} = \frac{\mm}{Kb^{(l)}} \text{ and }\,\, b^{(l+1)} = \frac{\mn}{Ka^{(l+1)}}.
\end{align}
This splitting, known as Sinkhorn's algorithm, is very efficient as it involves matrix vector multiplications only. Since the entropic regularization term \eqref{e:entropy} introduces blurring in the otherwise sparse solution, one is interested in keeping $\epsilon$ as small as  possible. Since the convergence of Sinkhorn's algorithm \eqref{e:sinkhorn} deteriorates as $\epsilon \rightarrow 0$, it is important to keep a balance between regularization and computational stability. In practice small values of $\epsilon$ lead to diverging scaling factors in \eqref{e:sinkhorn} and subsequent numerical instabilities. These problems can often be remedied using suitable scalings, see \cite{S2018}.

If the transportation costs depend on an underlying discrete structure, 
such as for our  graph-based migration problem, then the computational
burden of computing this cost must be take into consideration. For
our example the total transportation cost corresponds to the sum of 
edge weights when between vertices traversed on the shortest path. 
Note that the transportation costs are not necessarily the same in 
both directions since we consider directed graphs. We use Dijkstra’s 
algorithm to compute the shortest path from one node to all others 
in the graph, see \cite{D1959}. Dijkstra’s algorithm is based on continuous 
updates of the shortest distance to a starting point, and excludes 
longer distances in updates. It is the graph-based methodology that
underpins the  fast marching method to solve the eikonal equation
\cite{S1999}.

\subsection{MCMC and Optimization}
\label{ssec:MCMC}

We propose the use of Markov Chain Monte-Carlo (MCMC) methods to sample 
from the posterior distribution \eqref{e:posterior}. For the user
interested simply in optimization the algorithm we propose may be
viewed as a stochastic optimization method to reduce the model-data
misfit. MCMC methods
originated with the seminal paper \cite{metropolis1953equation} in which
what is now termed the  The Random walk Metropolis (RwM) algorithm was
introduced for a specific high dimensional integral required
in statistical physics. In our context the key desirable feature of
the method is that it requires only solution of the forward OT (or
regularized OT) problem, together with the generation of random numbers.
Given a current (approximate) sample from the posterior distribution, a
new sample is proposed by adding a mean zero Gaussian to the current one.
This is rejected if the resulting new state leaves $\sU$, and otherwise
accepted with a probability designed to preserve detailed balance with
respect to the posterior. The covariance of the Gaussian is an important 
tuning parameter: intuitively it should be chosen such that the acceptance rate
is neither close to $0$ or $1$, as either of these limits lead to 
successive iterates which are highly  correlated. The optimal scaling 
of RWM algorithms for different target densities has been investigated 
in \cite{RGG1997,RR2001}; although the theory developed there applies
in rather restricted scenarios, widespread experience and a variety
of theories demonstrate that
the work leads to useful rule-of-thumb for tuning acceptance
probabilities within the RwM algorithm \cite{YRR2019}, 
arguably because it leads to 
average acceptance probabilities that stay away from $0$ or $1$. 

In 1970 Hastings introduced a wide class of MCMC methods, now known
as Metropolis-Hastings algorithms \cite{hastings1970monte} and in
principle this provides a wide-range of variants on RwM that may be
used for our Bayesian formulation of inverse OT.
A popular variation of MCMC that we have found useful in the inverse
OT setting is Gibbs sampling. In high dimensional spaces it can be hard to
design proposals which are accepted with a reasonable acceptance probability,
and the idea of fixing subsets of the variables, and proposing moves in the
remainder, is natural. The Gibbs sampler allows this to be achieved in
a statistically consistent fashion. At each iteration one (or several) components of the unknown parameter is updated by sampling from its full conditional probability distribution, and cycling through all the variables. The method
may be relaxed to allow a RwM step from the conditional
probability distribution, rather than a full sample. The corresponding RwM-within-Gibbs method is outlined in Algorithm \ref{a:rwmwgibbs}. In this
algorithm we consecutively update $u$, $v$ and $W$ (or $f$). We generate proposals for each variable, which we accept or reject. Note that in general, for all 
the methods described here, any proposal which
descreases the value of $\Phi$ and remains in $\sU$ is accepted with 
probability one. Thus the Algorithm \ref{a:rwmwgibbs} may be viewed as an optimization
method which induces a stochastic gradient; the numerics will demonstrate
that this acts to minimize the misfit.

{\SetAlgoNoLine
  \begin{algorithm}
    \DontPrintSemicolon
  \caption{Random walk Metropolis within Gibbs}\label{a:rwmwgibbs}
  Initialize $(u^0, v^0, W^0)$\;
  \For{$k \geq 0$}{
    
    \KwSty{Generate} $\xi_u \sim \mathcal{N}(0,\delta_u^2)$ and \KwSty{propose} new value $x = u^k + \xi_u,~y = v^k, ~Z = W^k$\;
    \lIf{$(x,y,Z) \notin \sU$}{$(u^{k+1}, v^{k+1}, W^{k+1}) = (u^k, v^k, W^k)$}
    \Else{
      \begin{align*}\hspace*{-1em}
        (u^{k+1}, v^{k+1}, W^{k+1}) &=
        \begin{cases}
          (x,y,Z) &\hspace*{-1em}\text{ with probability } a((u^k, v^k,W^k),(x,y,Z))\\
           (u^k, v^k, W^k) &\hspace*{-1em}\text{ otherwise}
          \end{cases}
       \end{align*}}
    \KwSty{Generate} $\xi \sim \mathcal{N}(0,\delta_v^2)$ and \KwSty{propose} new value $y = v^k + \xi_v,~x = u^k, ~Z = W^k$\;
    \lIf{$(x,y,Z) \notin \sU$}{$(u^{k+1}, v^{k+1}, W^{k+1}) = (u^k, v^k, W^k)$}
    \Else{
      \begin{align*}
\hspace*{-1em}        (u^{k+1}, v^{k+1}, W^{k+1}) =
        \begin{cases}
          (x,y,Z) ~ &\hspace*{-1em}\text{with probability } a((u^k, v^k,W^k),(x,y,Z))\\
           (u^k, v^k, W^k) ~ &\hspace*{-1em}\text{otherwise}
          \end{cases}
       \end{align*}}
       \KwSty{Generate} $\xi_W \sim \mathcal{N}(0,\delta_W^2)$ and \KwSty{propose} new value $Z = W^k + \xi_W,~x^k = u^k, ~y = v^k$\;
    \lIf{$(x,y,Z) \notin \sU$}{$(u^{k+1}, v^{k+1}, W^{k+1}) = (u^k, v^k, W^k)$}
    \Else{
      \begin{align*}
      \hspace*{-1em}  (u^{k+1}, v^{k+1}, W^{k+1}) =
        \begin{cases}
          (x,y,Z) \quad &\hspace*{-1em}\text{with probability } a((u^k, v^k,W^k),(x,y,Z))\\
           (u^k, v^k, W^k) \quad &\hspace*{-1em}\text{otherwise}
          \end{cases}
       \end{align*}}
    Here
    \begin{align*}
      a((u,v,W),&(x,y,Z))=\\
      & \quad\min\Bigl\{1,\exp\bigl(\frac{1}{2\sigma^2}\lvert T - \cG(u,v,W)\rvert^2 -\frac{1}{2\sigma^2}\lvert T - \cG(x,y,Z)\rvert^2\bigr) \Bigr\}.
      \end{align*}
    }
\end{algorithm}
}

\section{Numerical Results}
\label{sec:NR}
In this section we demonstrate the behavior of MCMC methods for inverse 
OT, and Algorithm \ref{a:rwmwgibbs} in particular. 
We start by presenting results for the migration flow example introduced at the beginning and
use it as a 'proof-of-concept' for the proposed framework. We then continue with systematic numerical investigation to study
 the identifiability of the cost matrix in a variety of
scenarios, as well as discussing the behavior of the proposed methodology.
We focus on the three cost criteria discussed in Section \ref{ssec:GC}: Toeplitz cost \ref{i:toeplitz}, non-symmetric cost \ref{i:general} and graph-based cost \ref{i:graph}. We use the following functions implemented in the POT library \cite{pot} to solve the linear program \eqref{e:ot} as well as its regularized version \eqref{e:regot}:
\begin{itemize}
 \item \textit{emd} - this solver for linear programs is based on the respective network OT flow formulation of the problem and was introduced in \cite{BPPH2011}.
 \item \textit{sinkhorn} - implements the Sinkhorn-Knopp scaling algorithm to solve the regularized OT problem \eqref{e:regot} as proposed in \cite{Cuturi2013}.
\end{itemize}
We test the proposed methodologies using simulated data as well as real migration data. In making simulated data we compute the optimal
transportation maps $T$ for a given set of vectors $\mm$, $\mn$ and $f$ and
add i.i.d. Gaussian noise with mean $0$ and variance $\sigma^2$, see \eqref{e:noise}. Note that the resulting perturbed map $T^*$ may have negative entries and is not an element of $\cPnn$. Therefore we set all negative entries to zero and normalize it, to ensure that it is an admissible solution.

We illustrate the performance of the methodologies with plots of the running means and the respective posterior distributions. All posterior distributions are calculated after $500,000$ RwM iterations with a burn-in of $300,000$. The performed numerical experiments indicate that this number is sufficient for the convergence of MCMC. Note that we always plot the scaled vectors and matrices (unless stated otherwise). The penalty $\bar{C}$ in \eqref{e:pen} is set to $10$. Numerical simulations 
show that its absolute magnitude does not influence the posterior distributions significantly once above a certain level.\\

\subsection{European migration flows}

\begin{figure}[h!]
  \centering
  \includegraphics[width=0.95\textwidth]{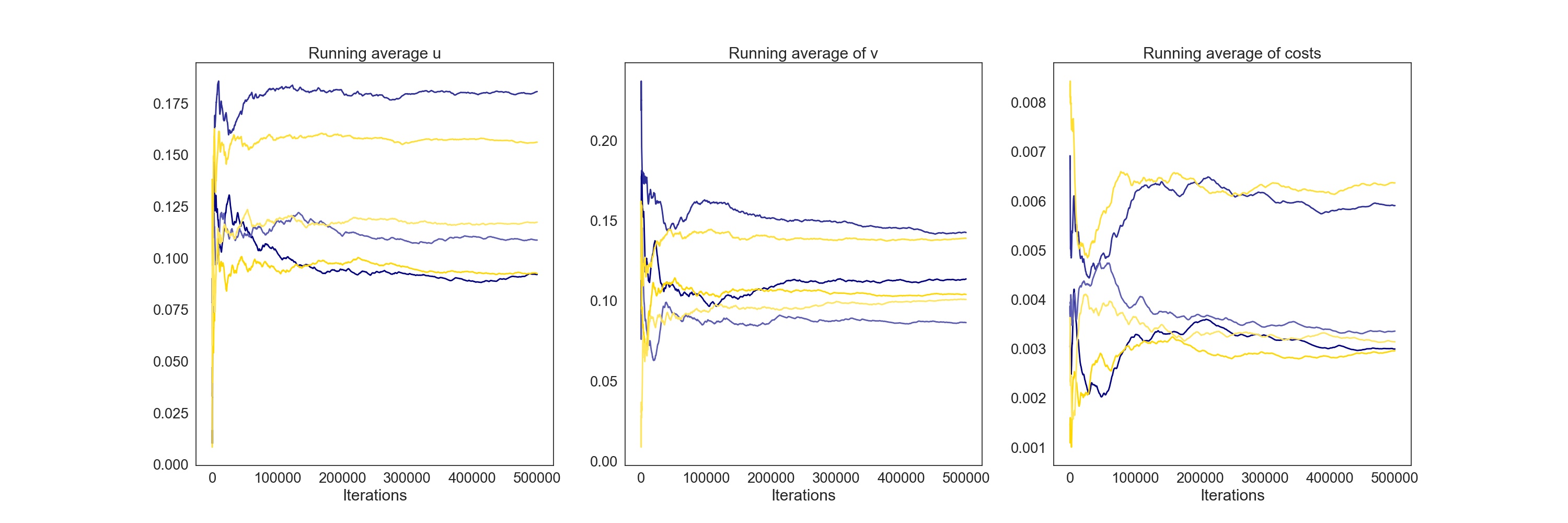}
  \caption{European network: Each plot shows the running average of three components of $u$, $v$ and $f$. The colors refer to the different combinations the exact LP solver (blue) and Sinkhorn (gold).}
  \label{f:eu_net_run_av}
\end{figure}

\begin{figure}[h!]
  \centering
  \begin{subfigure}[b]{0.3\textwidth}
      \includegraphics[width=0.95\textwidth]{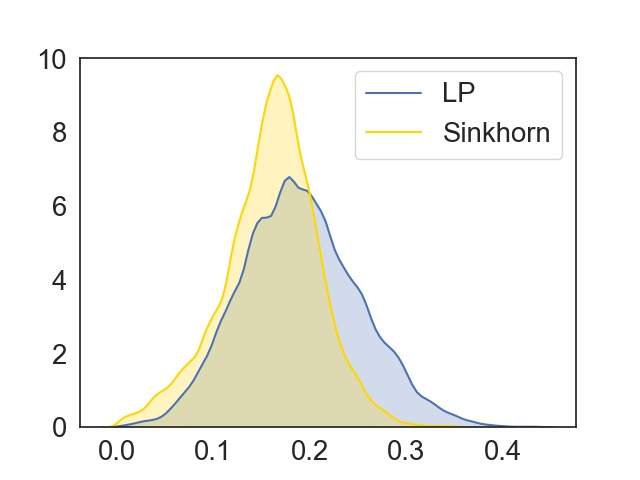}
      \caption{Posterior DE}
  \end{subfigure}
  \begin{subfigure}[b]{0.3\textwidth}
      \includegraphics[width=0.95\textwidth]{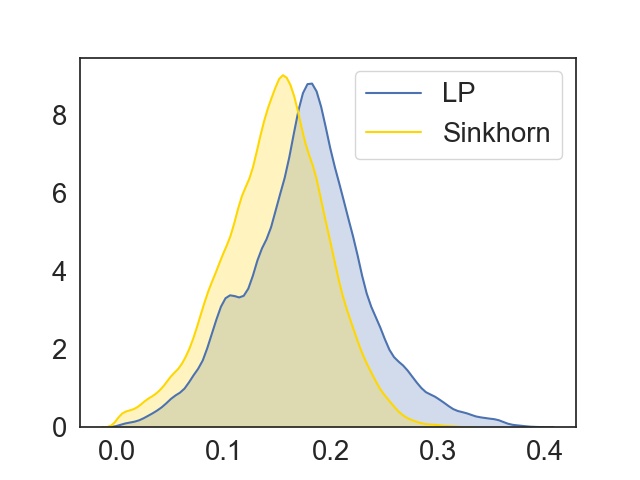}
      \caption{Posterior FR}
  \end{subfigure}
  \begin{subfigure}[b]{0.3\textwidth}
      \includegraphics[width=0.95\textwidth]{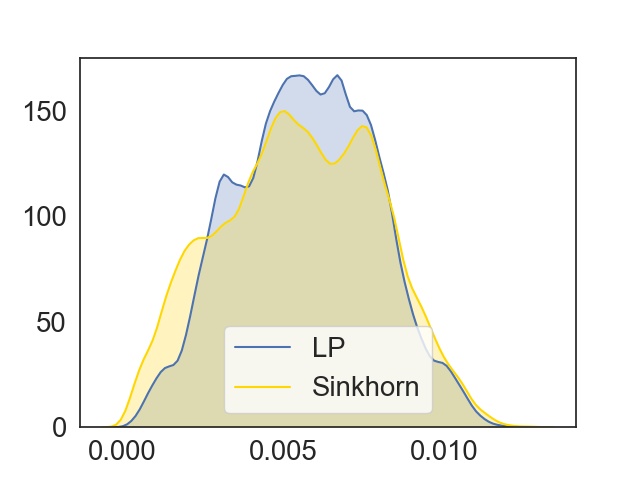}
      \caption{Posterior DE$\rightarrow$FR}
  \end{subfigure}
  \caption{European network: Posterior distributions of components of $u$, $v$ and $f$ using the exact LP solver and Sinkhorn.}
  \label{f:eu_net_posterior}
\end{figure}

\begin{figure}
  \centering
  \includegraphics[width=0.5\textwidth]{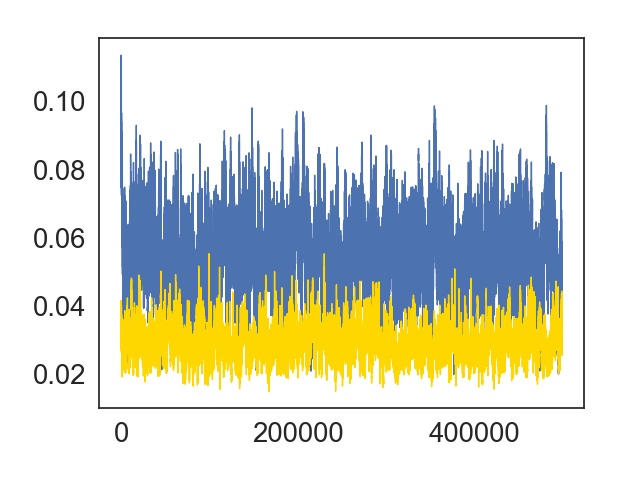}
  \caption{European network: misfit function for the exact LP solver (blue) and Sinkhorn (gold)}
  \label{f:eu_net_misfit}
  \end{figure}

We start by presenting estimates for the European network shown in Figure
\ref{f:europe}. We recall that vertices represent countries and that edges
connect countries sharing a border. The weights of these edges
correspond to the cost of moving from one country to another. The network
shown in Figure \ref{f:europe} consists of $n=9$ countries, which are
connected by $m=30$ edges. We use the estimated transportation map reported in
\cite{RWFSB2013} and assume that the noise level is $4\%$. The variance for
the proposals is set to $\delta_u^2 = \delta_v^2 = \delta_W ^2 =0.04$. We perform two
runs of the RwM-within-Gibbs algorithm, using the exact solver in the first and Sinkhorn's algorithm with $\epsilon=0.04$ in the second.
The acceptance rate of the exact solver is $50.8\%$ (($53.8\%$, $53.7\%$, $44.9\%$) for the different components $u$, $v$ and $f$), for Sinkhorn we have $82.9\%$ ($84.7\%$, $85.5\%$, $78.6\%$).
The running average of three components of  $u$, $v$ and $f$ are shown in Figure \ref{f:eu_net_run_av} and the corresponding posterior distributions in Figure \ref{f:eu_net_posterior}. We observe that both runs give comparable results, however the misfit for Sinkhorn is smaller, see Figure \ref{f:eu_net_misfit}. This 
difference might be explained by the fact that we underestimate the noise level $\sigma$ or that the actual transportation maps look more like solutions of 
regularized OT problems than the OT problem itself. 

\subsection{Graph-Based Cost}

Next we investigate the behavior of the proposed methodologies for graph-based
cost more thoroughly. We will see that
\begin{itemize}
\item The identification  of $u$, $v$ and $f$ is robust with respect to the
  sampling variances, see Figure \ref{f:toy_graph_posterior}.
\item The posterior estimates are consistent using different solvers, see
  Figure \ref{f:toy_graph_posterior_2}.
  \end{itemize}
These results are obtained using noisy transportation maps $T^*$ for a graph connecting $n=5$ nodes with $m=12$ edges. In doing so we solve problem \eqref{e:ot}
for given vectors $\mm$, $\mn$ and $f$ and add $4\%$ noise. Note that this inverse problem is
overdetermined since $2 \cdot 5 + 12 - 3 < 5^2 -1$.

\subsubsection*{Influence of the Sampling Variance $\delta^2$} 
\begin{figure}[h!]
  \centering
  \includegraphics[width=0.95\textwidth]{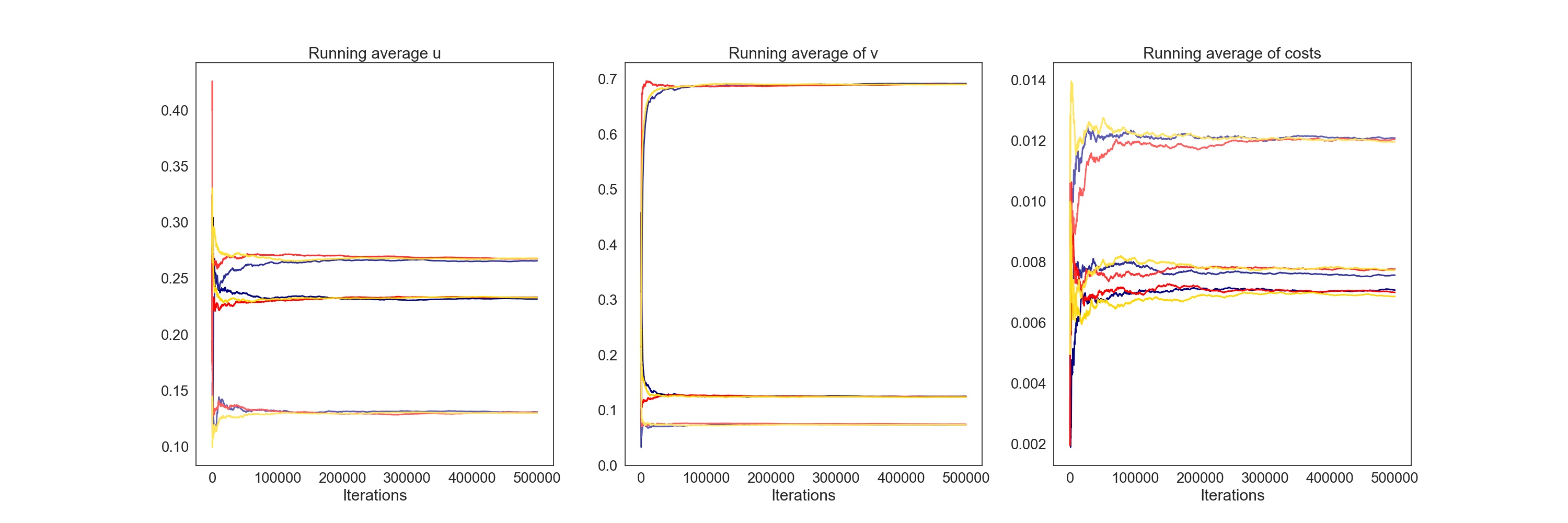}
  \caption{Graph based cost: Each plot shows the running average of three components of $u$, $v$ and $f$. The colors refer to the different combinations of $\delta$ used - blue $(\delta_u^2, \delta_v^2, \delta_f^2) = (0.02,0.02,0.04)$, red $(\delta_u^2, \delta_v^2, \delta_f^2) = (0.04,0.04,0.04)$ and gold $(\delta_u^2, \delta_v^2, \delta_f^2) = (0.04,0.02,0.04)$.}
  \label{f:toy_graph_run_av}
\end{figure}

\begin{figure}[h!]
  \centering
  \begin{subfigure}[b]{0.3\textwidth}
      \includegraphics[width=0.95\textwidth]{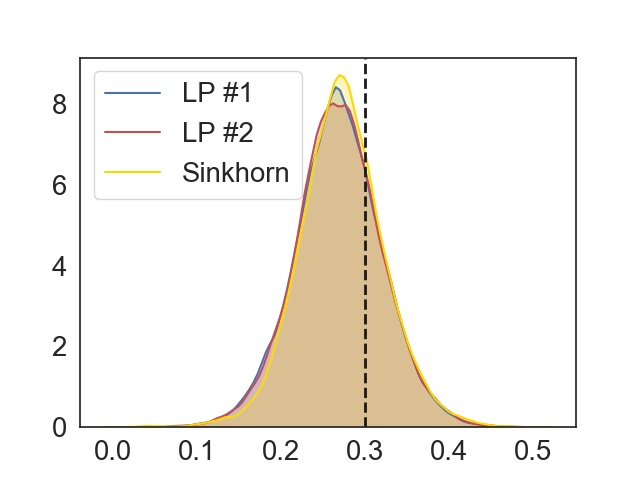}
      \caption{Posterior $u_2$}
  \end{subfigure}
  \begin{subfigure}[b]{0.3\textwidth}
      \includegraphics[width=0.95\textwidth]{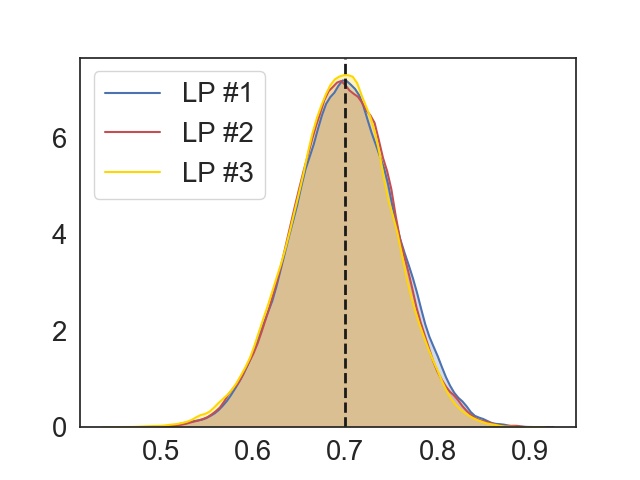}
      \caption{Posterior $v_2$}
  \end{subfigure}
  \begin{subfigure}[b]{0.3\textwidth}
      \includegraphics[width=0.95\textwidth]{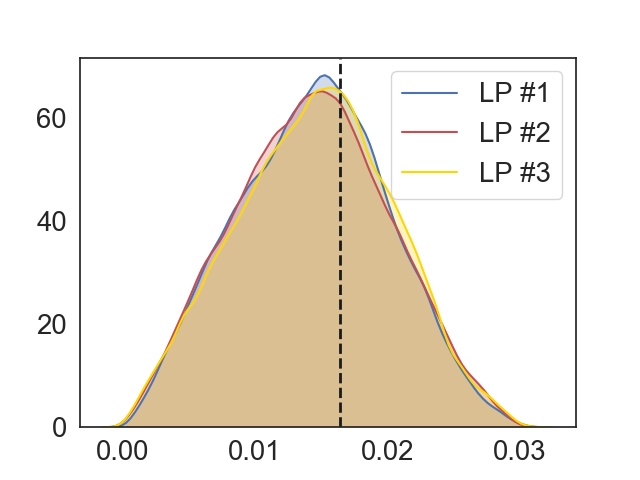}
      \caption{Posterior $C_{1,2}$}
  \end{subfigure}
  \caption{Graph based cost: Posterior distributions of components of $u$, $v$ and $f$ using different combinations of $\delta$.}
  \label{f:toy_graph_posterior}
\end{figure}

We start by
investigating the impact of the sampling variance $\delta^2$. We perform MCMC
runs for different combinations of $\delta_u$, $\delta_v$ and $\delta_f$ (listed in Table \ref{t:graph_delta_accept}) and
compute the running average and posterior distributions of some
components. Note that these parameters
affect the rate of convergence of the algorithm, but not the posterior
distribution itself. The variance of the samples determines how much new samples differ from the previous iterates - the larger the variance the more 
adventurous the search, but the less likely to accept leading to highly
correlated samples because of rejections. On the other hand
smaller variance has a higher probability of accepting but is not adventurous
and hence leads to highly correlated samples. It is thus desirable to have an
acceptance rate that is neither close to $0$ or $1$. The running averages of three components of $u$, $v$ and $W$ are
shown in Figure \ref{f:toy_graph_run_av}, the respective posteriors in
Figure \ref{f:toy_graph_posterior}. We see that the results are consistent
for all combinations of $\delta$'s. However the respective convergence rates vary, see
Table \ref{t:graph_delta_accept}. We observe a generally higher acceptance rate when
sampling from the marginal distribution of $\mm$, and a decreased rate when
increasing the sampling variance.

\begin{table}
  \centering
  \begin{tabular}{*{3}{|c}||*{4}{c|}}
    \hline
  $\delta_u^2$ & $\delta_v^2$ & $\delta_f^2$ & $a$ & $a_u$ & $a_v$ & $a_f$  \\
  \hline
  0.02 & 0.02 & 0.04 & 65 & 80 & 52 & 62 \\
    0.04 & 0.04 & 0.04 & 51 & 65 & 26 & 62 \\
    0.04 & 0.02 & 0.04 & 60 & 65 & 52 & 62\\
    \hline
\end{tabular}
\caption{Graph based cost: acceptance rates  in $\%$ for different combinations of $\delta_u$, $\delta_v$ and $\delta_f$.}
\label{t:graph_delta_accept}
\end{table}

\subsubsection*{Exact vs. Sinkhorn} 
\begin{figure}[h!]
  \centering
  \begin{subfigure}[b]{0.3\textwidth}
      \includegraphics[width=0.95\textwidth]{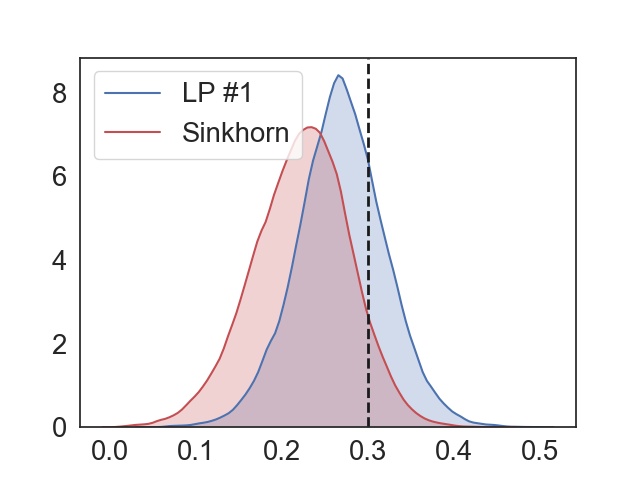}
      \caption{Posterior $u_2$}
  \end{subfigure}
  \begin{subfigure}[b]{0.3\textwidth}
      \includegraphics[width=0.95\textwidth]{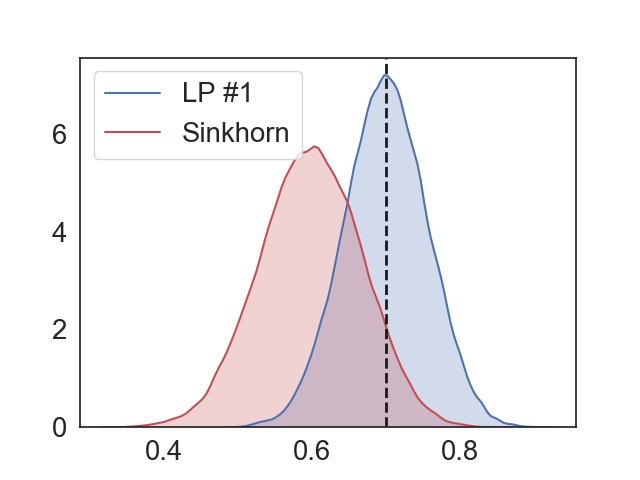}
      \caption{Posterior $v_2$}
  \end{subfigure}
  \begin{subfigure}[b]{0.3\textwidth}
      \includegraphics[width=0.95\textwidth]{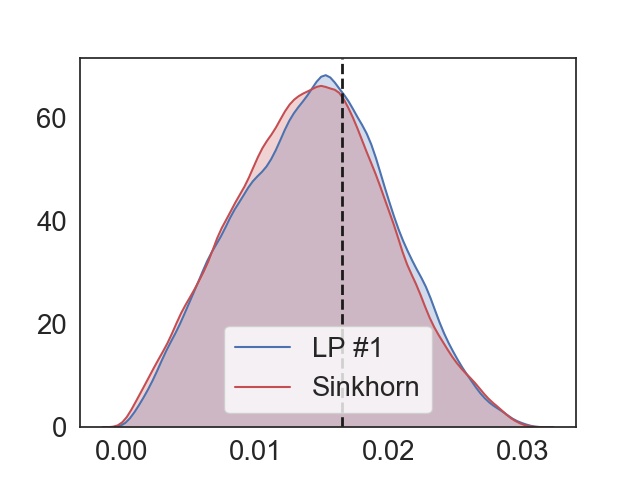}
      \caption{Posterior $C_{1,2}$}
  \end{subfigure}
  \caption{Graph based cost: Posterior distributions of components of $u$, $v$
    and $f$ when using the exact LP solver in the RwM with Gibbs, or regularizing
the LP to employ the Sinkhorn algorithm.}
  \label{f:toy_graph_posterior_2}
\end{figure}

Next we investigate the sensitivity of
the results with respect to the forward solver used in Algorithm
\ref{a:rwmwgibbs}. We run two RwM simulations - the first one using the
exact solver and the second one using the Sinkhorn algorithm. We observe that
both runs give similar posterior distributions if we choose the regularization
parameter $\epsilon$ in a sensible way, see Figure
\ref{f:toy_graph_posterior_2}. Generally speaking it seems advisable
to choose it similar to the noise level (as in the shown results). We will investigate the impact of the
regularization parameter in the next subsection in more detail.
 
\subsection{Toeplitz Cost}

In the following we present more detailed numerical experiments if $C$ is Toeplitz. The findings 
of the numerical experiments performed in this subsection can be summarized as follows:
\begin{itemize}
\item The posterior distributions of $u$, $v$
  and $f$ are consistent for varying ranges of proposal variances $\delta$, see
  Figure \ref{f:toeplitz_run_av_var_sig} and Figure \ref{f:toeplitz_posterior_var_sig}.
  \item The exact solver and Sinkhorn's algorithm converge to similar
    posteriors, if the entropic regularization parameter $\epsilon$ is chosen
    sensibly, see Figure \ref{f:toeplitz_run_av} and Figure
    \ref{f:toeplitz_posterior}.
  \item The variance of the posteriors increases with the noise level in the
    data, as shown for example in Figure \ref{f:toeplitz_posterior_2} and
    Figure \ref{f:toeplitz_posterior_3}.
    \item Sinkhorn's algorithm gives a higher acceptance rate and a more monotone decrease of the
     data-misfit function, see  Figure \ref{f:toep_sink_vs_emd_misfit}.
\end{itemize}
We underpin these statements with numerical simulations using generated noisy
transportation maps.  We recall that $C$ has $2n-3$ degrees of
freedom in case of Toeplitz cost \ref{i:toeplitz}. This defines, as in the
case of graph-based cost \ref{i:graph}, a mapping from the 
vector $f \in \mathbb{R}^{2n-3}$ to the cost matrix $C$, that is $\mathcal{E}: \mathbb{R}^{2n-3} \rightarrow
\cPnn$ with $C = \mathcal{E}(f)$. Hence we
generate proposals for the vector $f$, which define the entries of $C$.\\ 
We set $n=5$ and generate a noisy realization $T^*$ for a given set of vectors $\mm$, $\mn \in
 \mathcal{P}_5$ and $f \in \mathbb{R}^7$ (which is mapped to the
 respective Toeplitz cost matrix $C$ $\mathcal{P}_{5 \times 5}$). Then $T^*$ is
 obtained by adding noise $\eta$ with variance $\sigma^2 = 0.04$ (unless
 stated otherwise) and subsequent normalization of the distorted map. Note
 that this problem is overdetermined, since $C$ is Toeplitz and $n > 2$.

\subsubsection*{Influence of the Sample Variance $\delta^2$} 
\begin{figure}[ht!]
  \centering
  \includegraphics[width=\textwidth]{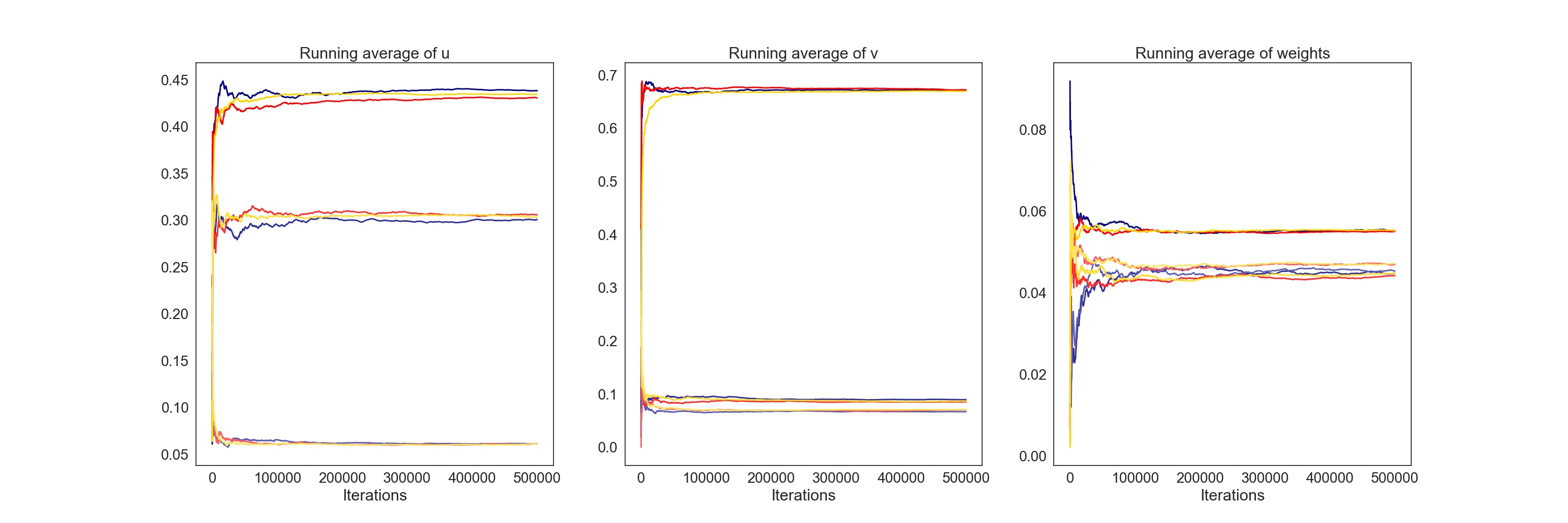}
  \caption{Toeplitz cost: Running averages of three components of $u$, $v$ and
    $f$. The colors refer to different variances of the proposals. The blue plots correspond to $\delta_u^2 = \delta_v^2 = \delta_f^2 = \delta^2 = 0.02$, the red ones to $\delta_u^2 = \delta_v^2 = 0.02$ and $\delta_f^2 = 0.04$ and the golden ones to $\delta_u^2 = \delta_v^2 = \delta_f^2 = 0.04$. }
  \label{f:toeplitz_run_av_var_sig}
\end{figure}

\begin{figure}[ht!]
  \centering
  \begin{subfigure}[b]{0.45\textwidth}
      \includegraphics[width=0.95\textwidth]{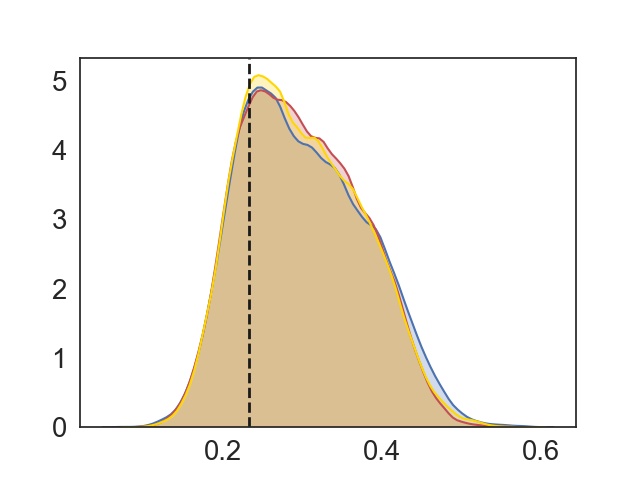}
      \caption{Posterior $u_3$}
  \end{subfigure}
  \begin{subfigure}[b]{0.45\textwidth}
      \includegraphics[width=0.95\textwidth]{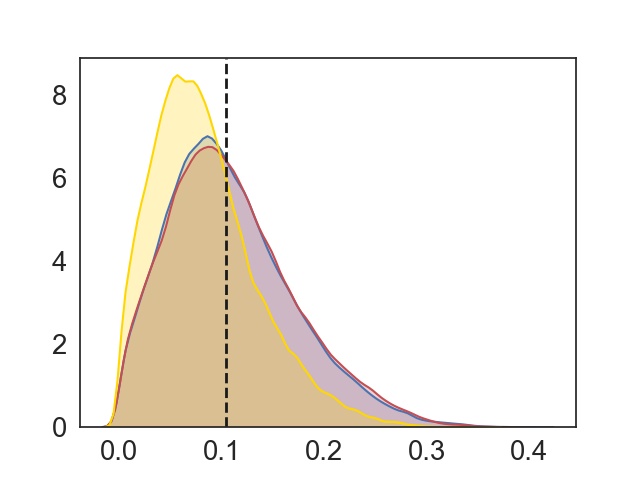}
      \caption{Posterior $v_1$}
  \end{subfigure}
  \begin{subfigure}[b]{0.45\textwidth}
      \includegraphics[width=0.95\textwidth]{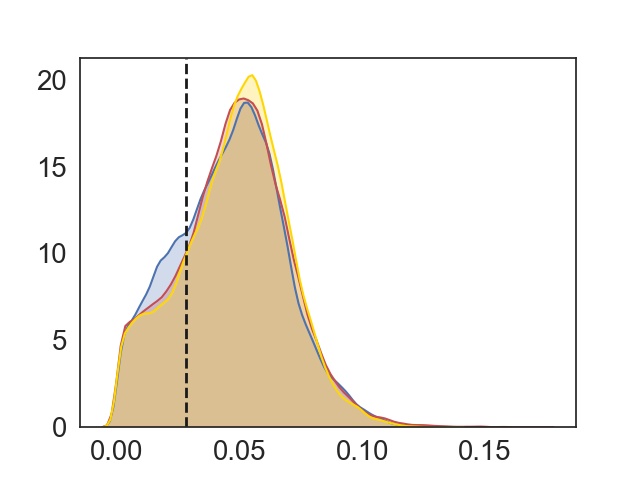}
      \caption{Posterior $f_4$}
  \end{subfigure}
  \begin{subfigure}[b]{0.45\textwidth}
      \includegraphics[width=0.95\textwidth]{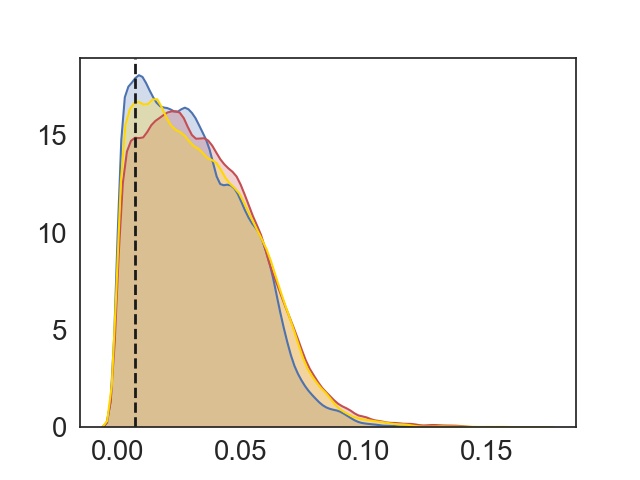}
      \caption{Posterior $f_8$}
  \end{subfigure}

  \caption{Toeplitz cost: posterior distributions of four components of $u$, $v$ and $f$. The colors refer to different variances of the proposals. 
The blue plots correspond to $\delta_u^2 = \delta_v^2 = \delta_f^2 = \delta^2 = 0.02$, the red ones to $\delta_u^2 = \delta_v^2 = 0.02$ and $\delta_f^2 = 0.04$ and the golden ones to $\delta_u^2 = \delta_v^2 = \delta_f^2 = 0.04$. }
  \label{f:toeplitz_posterior_var_sig}
\end{figure}

As in the case of graph-based cost, we investigate  the performance of the RwM-within-Gibbs algorithm for different
combinations of  $\delta_u$, $\delta_v$ and $\delta_W$.  Table
\ref{t:toeplitz_delta_accept} lists the considered $\delta$-combinations
together with the acceptance rates. The running average of the $3$ or $4$
different components of the posteriors are shown in Figures
\ref{f:toeplitz_run_av_var_sig}  and \ref{f:toeplitz_posterior_var_sig}. The figures show, as expected, that the posterior
distributions are independent of the choice of the $\delta$ parameters.
They also show that, in the ranges chosen, the rate of convergence does
not vary in any considerable way -- the method is fairly robust.

\begin{table}
  \centering
  \begin{tabular}{*{3}{|c}||*{4}{c|}}
    \hline
  $\delta_u^2$ & $\delta_v^2$ & $\delta_f^2$ & $a$ & $a_u$ & $a_v$ & $a_f$  \\
  \hline
  0.02 & 0.02 & 0.02 & 66.1 & 67.9 & 55.4 & 75 \\
    0.02 & 0.02 & 0.04 & 60.3 & 68.3 & 55.3 & 52.7 \\
    0.04 & 0.04 & 0.04 & 44.1 & 45.5 & 30.0 & 57\\
    \hline
\end{tabular}
\caption{Toeplitz cost: acceptance rates for different combinations of $\delta_u$, $\delta_v$ and $\delta_f$.}
\label{t:toeplitz_delta_accept}
\end{table}

\subsubsection*{Exact vs. Sinkhorn:} 
\begin{figure}[ht!]
  \centering
  \includegraphics[width=\textwidth]{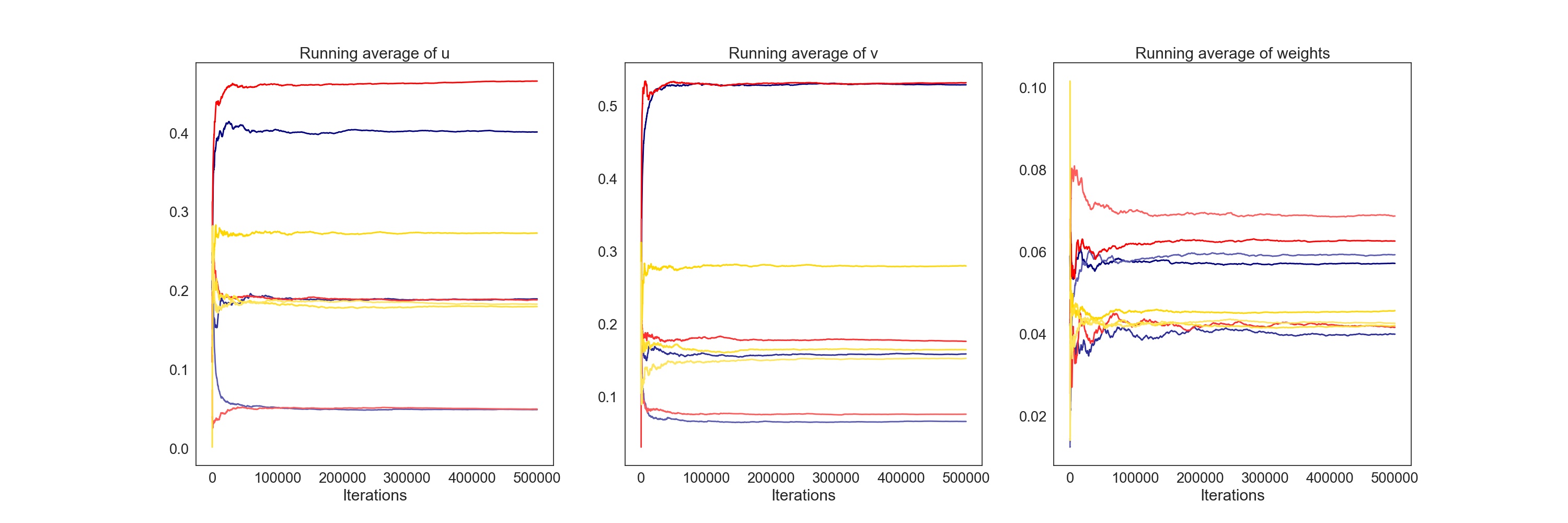}
  \caption{Toeplitz cost: Running averages of three components of $u$, $v$ and $f$. The colors refer to the used solver for the forward OT problem. Red corresponds to the exact LP solver, blue and gold to the Sinkhorn algorithm with regularization parameter $\epsilon = 0.04$ and $\epsilon=0.1$ respectively.}
  \label{f:toeplitz_run_av}
\end{figure}

\begin{figure}[ht!]
  \centering
  \begin{subfigure}[b]{0.32\textwidth}
      \includegraphics[width=0.95\textwidth]{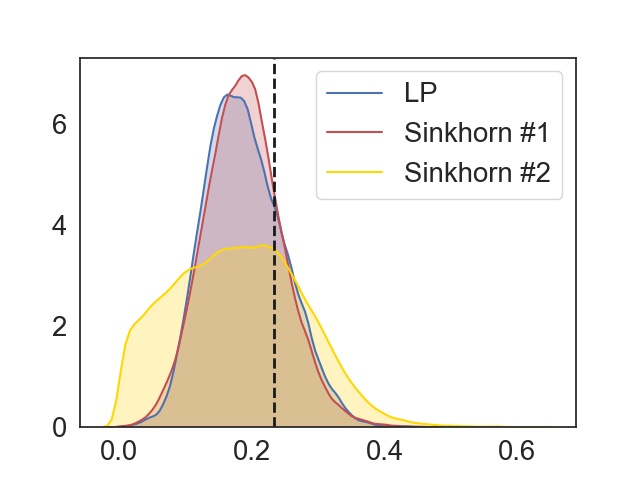}
      \caption{Posterior $u_3$}
  \end{subfigure}
  \begin{subfigure}[b]{0.32\textwidth}
      \includegraphics[width=0.95\textwidth]{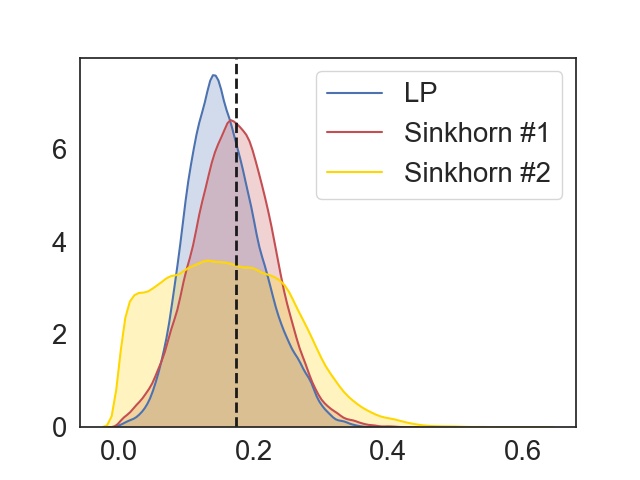}
      \caption{Posterior $v_3$}
  \end{subfigure}
  \begin{subfigure}[b]{0.32\textwidth}
      \includegraphics[width=0.95\textwidth]{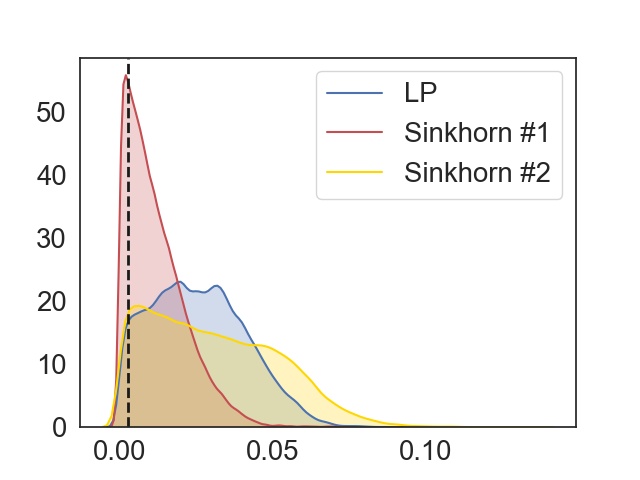}
      \caption{Posterior $C_{1,2}$}
  \end{subfigure}
  \caption{Toeplitz cost: posterior distributions of three components of $u$, $v$ and $f$. The colors refer to the solver used for the forward OT problem. Red corresponds to the exact LP solver, blue and gold to the Sinkhorn algorithm with regularization parameter $\epsilon = 0.04$ and $\epsilon=0.1$ respectively.  }
  \label{f:toeplitz_posterior}
\end{figure}

Next we take a closer look how results change if
we use Sinkhorn's algorithm instead of the exact solver. In particular we
investigate how the size of the regularization parameter $\epsilon$ as well as
the way we generate data effects the performance and results of the RwM algorithm.\\
We start by generating a noisy transportation map using the exact solver for
\eqref{e:ot}. Then we compare the posterior distributions using the exact
solver for the reconstruction in the first run and  the Sinkhorn algorithm with $\epsilon = 0.04$ and $\epsilon=0.1$ in the next two test runs. 
Figure \ref{f:toeplitz_run_av} shows the running average of three components
of the vectors $\mm$, $\mn$ and $f$ (left to right). The color coding relates
to the solver used - red corresponds to the exact forward solver, blue and
yellow when the Sinkhorn algorithm was used. Figure \ref{f:toeplitz_posterior}
shows the posterior distribution of the second component of $\mm$ and $\mn$ as
well as the fifth entry of the vector $f$. We observe that we obtain similar
posteriors when using the exact solver (LP) and Sinkhorn with $\epsilon =
0.04$. If the regularization parameter $\epsilon$ is chosen larger, which
results in blurred  (and therefore less sparse) transportation maps, the posterior
distributions are less pronounced and close to uniform on the respective
scaled intervals (due to the normalization constraint).

\begin{figure}[ht!]
  \centering
  \begin{subfigure}[b]{0.32\textwidth}
      \includegraphics[width=0.95\textwidth]{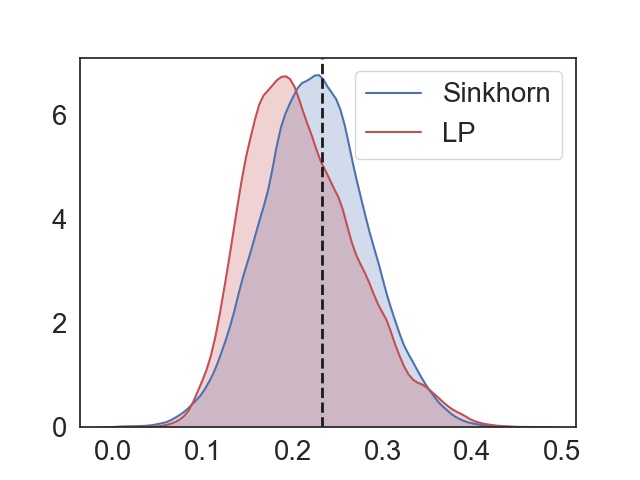}
      \caption{Posterior $u_3$}
  \end{subfigure}
  \begin{subfigure}[b]{0.32\textwidth}
      \includegraphics[width=0.95\textwidth]{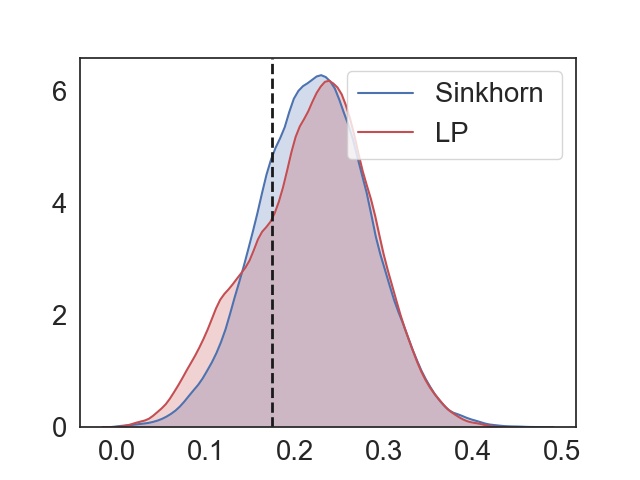}
      \caption{Posterior $v_3$}
  \end{subfigure}
  \begin{subfigure}[b]{0.32\textwidth}
      \includegraphics[width=0.95\textwidth]{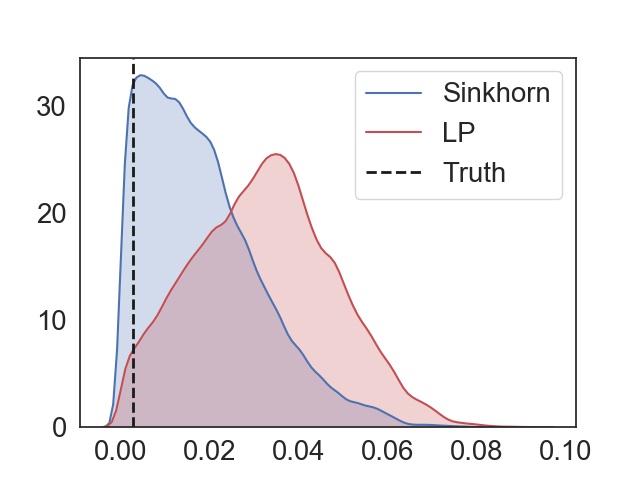}
      \caption{Posterior $C_{1,2}$}
  \end{subfigure}
  \caption{Toeplitz cost: Posteriors of $u_3$, $v_3$ and $C_{1,2}$ using data generated by the Sinkhorn algorithm with $\epsilon = 0.04$.}
  \label{f:toeplitz_posterior_2}
\end{figure}

\begin{figure}[ht!]
  \centering
  \begin{subfigure}[b]{0.32\textwidth}
      \includegraphics[width=0.95\textwidth]{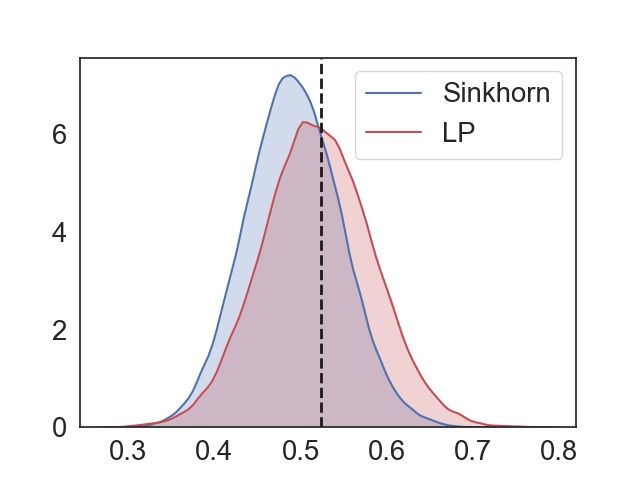}
      \caption{Posterior $u_1$}
  \end{subfigure}
  \begin{subfigure}[b]{0.32\textwidth}
      \includegraphics[width=0.95\textwidth]{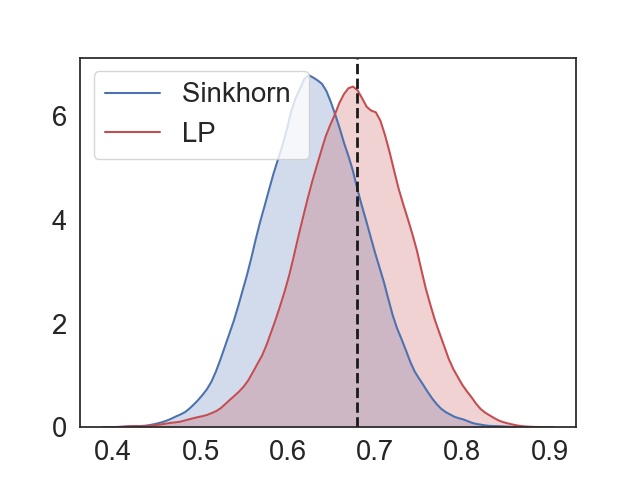}
      \caption{Posterior $v_2$}
  \end{subfigure}
  \begin{subfigure}[b]{0.32\textwidth}
      \includegraphics[width=0.95\textwidth]{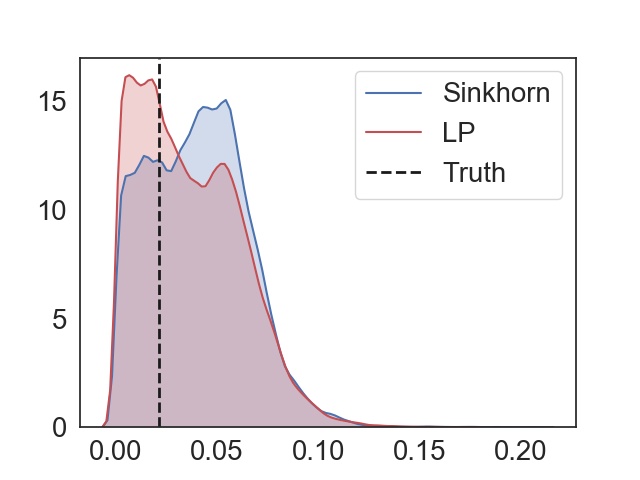}
      \caption{Posterior $C_{5,1}$}
  \end{subfigure}
  \caption{Toeplitz cost: Posteriors of $u_1$, $v_2$ and $C_{5,1}$ using data generated by the Sinkhorn algorithm with $\epsilon = 0.1$.}
  \label{f:toeplitz_posterior_3}
\end{figure}

\begin{figure}
  \centering
  \begin{subfigure}[b]{0.475\textwidth}
      \includegraphics[width=0.95\textwidth]{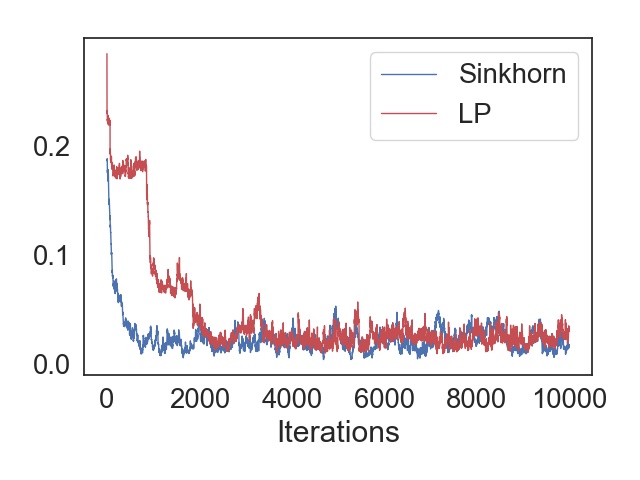}
  \caption{Data generated with $\epsilon=0.04$.}
  \end{subfigure}
  \begin{subfigure}[b]{0.475\textwidth}
      \includegraphics[width=0.95\textwidth]{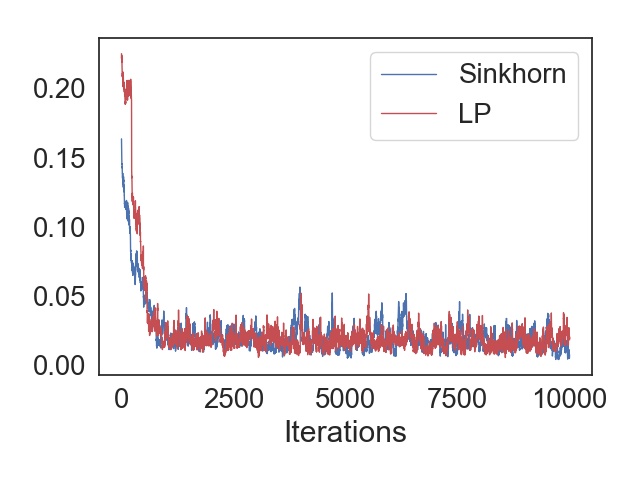}
  \caption{Data generated with $\epsilon=0.1$.}
  \end{subfigure}
  \caption{Toeplitz cost: First $10~000$ iterations of the misfit function \eqref{e:misfit}.}
  \label{f:toep_sink_vs_emd_misfit}
\end{figure}

Next we generate the noisy transportation map using the Sinkhorn algorithm. We
set the regularization parameter $\epsilon=0.04$ and we distort the computed
map with $4\%$ and $10\%$ noise. In each case we perform two different RwM runs, first using the Sinkhorn
algorithm and then the exact solver. The respective posterior distributions
are shown in Figure \ref{f:toeplitz_posterior_2} and Figure
\ref{f:toeplitz_posterior_3}.
We observe no significant difference in the quality of the posteriors. Figure
\ref{f:toep_sink_vs_emd_misfit} illustrates an interesting difference in the
convergence behavior of the RwM algorithm. The data misfit term
\eqref{e:misfit} shows multiple drops when using the exact
solver. Such jumps haven not been observed when using the Sinkhorn
algorithm. We recall that the Sinkhorn algorithm solves the
respective regularized (convex) optimization problem, which has a unique
minimum. We believe that the non-uniqueness of the exact forward problem leads
to several local minima in the inverse problem, in which the RwM algorithm gets stuck.

\begin{figure}[ht!]
  \centering
  \includegraphics[width=\textwidth]{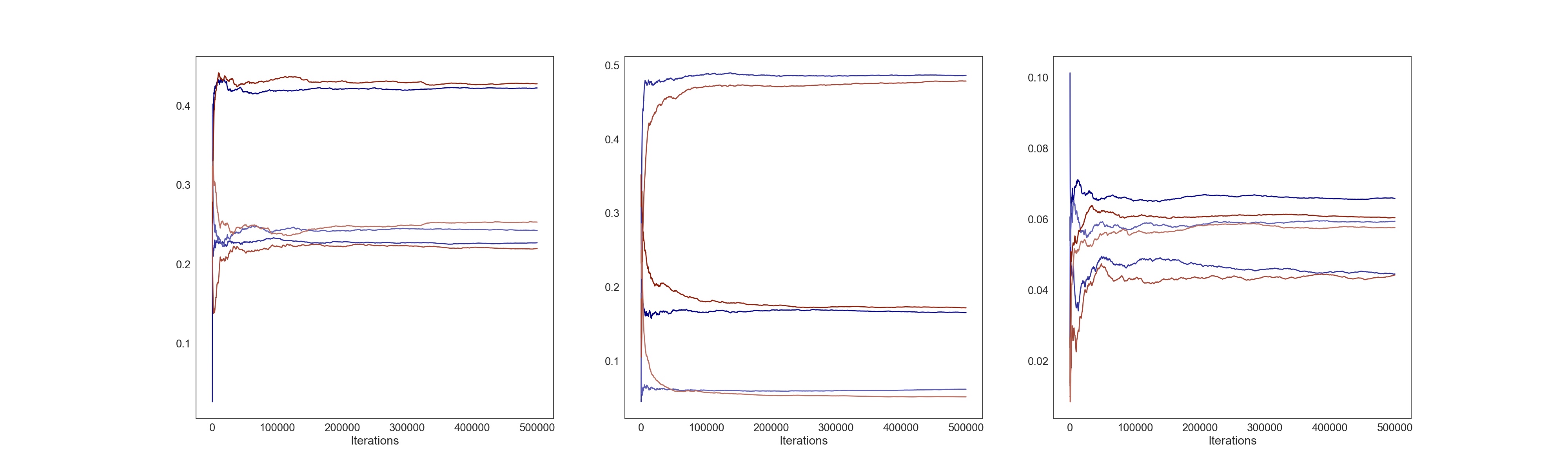}
  \caption{Toeplitz cost: Running averages of three components of $u$, $v$ and $f$ for data generated by the Sinkhorn algorithm with $\epsilon=0.04$. The blue plots are the running averages using Sinkhorn in the RwM, the red ones the exact LP solver. }
  \label{f:toeplitz_run_av2}
\end{figure}

\subsection{General Cost}
\begin{figure}[ht!]
  \centering
  \begin{subfigure}[b]{0.32\textwidth}
      \includegraphics[width=0.95\textwidth]{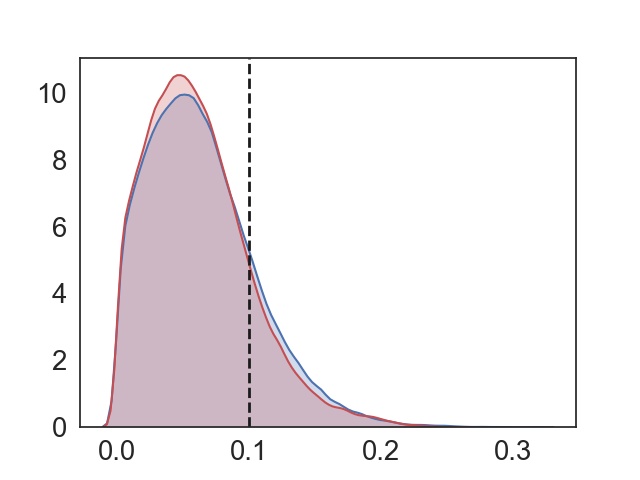}
      \caption{Posterior $u_2$}
  \end{subfigure}
  \begin{subfigure}[b]{0.32\textwidth}
      \includegraphics[width=0.95\textwidth]{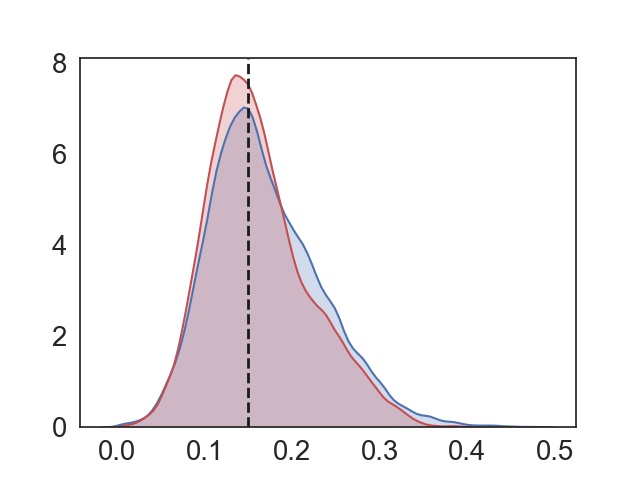}
      \caption{Posterior $v_2$}
  \end{subfigure}
  \begin{subfigure}[b]{0.32\textwidth}
      \includegraphics[width=0.95\textwidth]{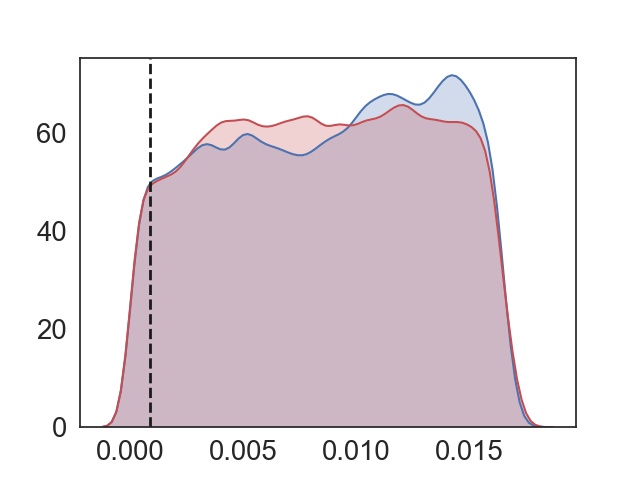}
      \caption{Posterior $C_{1,2}$}
  \end{subfigure}
  \caption{General cost: posterior distributions of different components of $u$, $v$ and $W$. The two colors refer to the different combinations of $\delta$ - red to  $\delta_u^2 = \delta_v^2 = \delta_W^2 = \delta = 0.02$ and ones to $\delta_u^2 = \delta_v^2 = 0.02$ and $\delta_W^2 = 0.04$.}
  \label{f:general_posterior}
\end{figure}

So far we investigated overdetermined problems only. Hence we conclude by
considering general non-symmetric costs, that is case \ref{i:general}, for
$n=5$. This identification problem is underdetermined and we expect poorer
identifiability and quality of posteriors. This presumption is confirmed by
our numerical experiments, see for example Figure \ref{f:general_posterior}.\\

\begin{figure}[ht!]
  \centering
  \includegraphics[width=\textwidth]{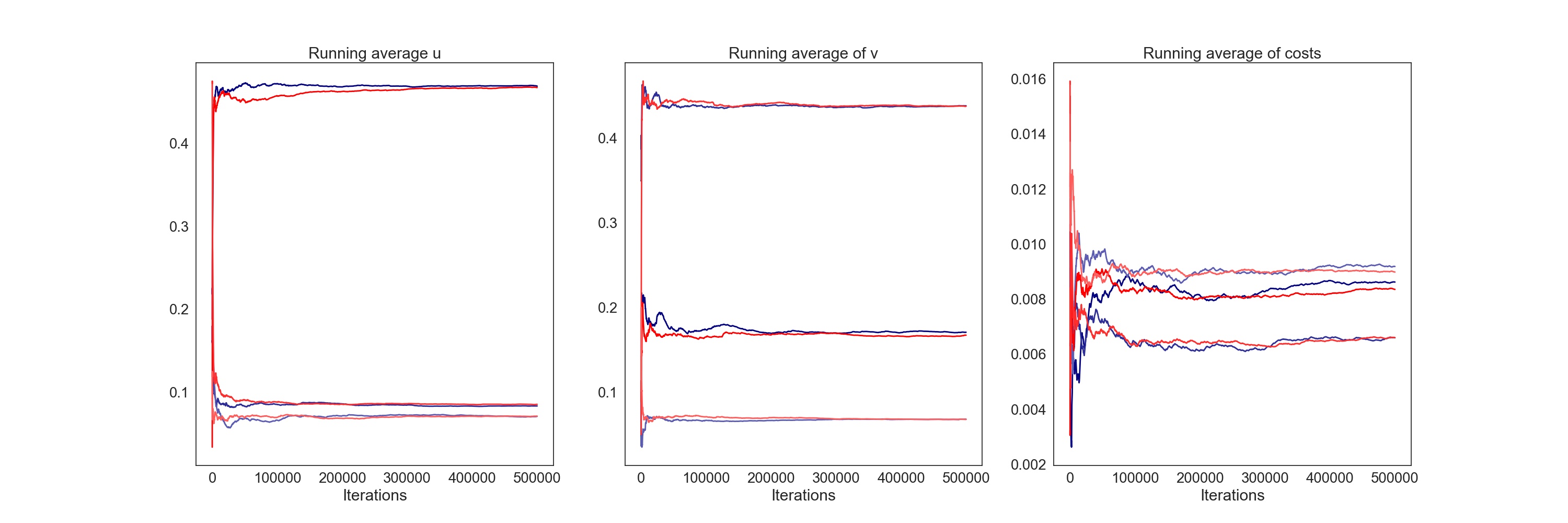}
  \caption{General cost: Each plot shows the running average of three components of $u$, $v$ and $W$. The colors correspond to different combinations of $\delta$ - red plots to $\delta_u^2 = \delta_v^2 = \delta_W^2 = \delta^2 = 0.02$ and blue ones to $\delta_u^2 = \delta_v^2 = 0.02$ and $\delta_W^2 = 0.04$.}
  \label{f:general_run_av}
\end{figure}

We investigate the identification from generated data in case of $4\%$ noise. We
perform two RwM test runs using the exact solver to calculate the posterior
distributions of $u$, $v$ and $W$. In the first run we set the sample variance to $\delta_u^2 =
\delta_v^2 = \delta_W^2 = \delta = 0.02$  and in the second  to $\delta_u^2 = \delta_v^2 = 0.02$ and $\delta_W^2 = 0.04$.
Figure \ref{f:general_run_av} shows the running averages for $3$ different
components of $u$, $v$ and $W$ for both runs. We observe that the components
of $u$ and $v$ converge much faster than the ones of $W$ and that the convergence
is consistent for both sets of $\delta$'s. The posterior
distributions of $u$ and $v$ give reasonable results, while the posteriors of
the cost matrix are close to uniform on the respective scaled intervals (due
to the normalization constraint). This indicates that the components of
the cost matrix $W$ are difficult to identify. We expect that the
identifiability gets worse as the dimension $n$ increases.

\section{Conclusions}
\label{sec:C}

This paper introduces a systematic approach to infer unknown costs from noisy
observations of optimal transportation plans. It is based on the Bayesian
framework for inverse problems and allows to  quantify uncertainty in the obtained
estimates; however the methodology may also be viewed as a stochastic
optimization procedure in its own right, tuning the unknowns
so that the optimal transport plan better fits the data. The performance
of the developed methodologies is investigated using the example of
international migration flows. In this context reported annual migration flow
statistics can be interpreted as noisy observations of optimal transportation
plans with cost related to the geographical position of countries. We
formulate the graph-based problem, estimate the weights, which
represent the costs of moving between neighbouring countries, and 
quantify uncertainty in the weights.
Our numerical investigation show that the proposed methodologies are robust and consistent for different cost
functions and parametrizations. We observed that the distributions as well as the costs can be accurately
determined for a variety of settings, if the
problem is overdetermined. The identifiability declines as
the dimensionality increases or if the problem becomes underdetermined.\\

The proposed framework provides the basis for a multitude of future
research directions in applied mathematics and other scientific
disciplines. The next steps will focus on several questions related to the use  of the Sinkhorn algorithm in the context of inverse optimal transport, such as the convergence rate of the regularized problem \eqref{e:regot} as $\epsilon \rightarrow 0$
or the optimal choice of $\epsilon$ with respect to the noise level $\sigma$;
furthermore hierarchical algorithms which learn parameters such as these
from the data would also be of interest.
In the context of migration flows, different modeling aspects, such
as the coupling to age structured population models or the formulation of the
OT problem on the continuous level, will be investigated. Furthermore the
application of the developed methodologies for general linear programs, which play an important role in transportation research, manufacturing,
economics and demography, will be of interest. \vspace*{1em}\\

\noindent{\textbf{Acknowledgments.} The authors are
grateful to Venkat Chandrasekaran
helpful discussions about the literature in inverse linear programming.
The work of AMS is funded by US National Science Foundation (NSF) grant DMS
1818977 and AFOSR Grant FA9550-17-1-0185. The work of MTW was partly supported
by The Royal Society International Exchanges grant IE 161662.}

\bibliographystyle{plain}
\bibliography{migration}


\end{document}